\documentclass[11pt,letterpaper,reqno]{amsart}
\usepackage[margin=1.2in]{geometry}

\usepackage[T1]{fontenc}
\usepackage{amsmath,amssymb}
\usepackage[mathscr]{euscript}
\usepackage{microtype}
\usepackage{tikz-cd}
\usepackage[dvipsnames]{xcolor}
\usepackage{aliascnt}
\usepackage{dsfont} 

\theoremstyle{plain}
\newtheorem{theorem}{Theorem}[section]
\newaliascnt{proposition}{theorem}
\newtheorem{proposition}[proposition]{Proposition}
\aliascntresetthe{proposition}
\newaliascnt{lemma}{theorem}
\newtheorem{lemma}[lemma]{Lemma}
\aliascntresetthe{lemma}
\newaliascnt{corollary}{theorem}
\newtheorem{corollary}[corollary]{Corollary}
\aliascntresetthe{corollary}
\newaliascnt{conjecture}{theorem}
\newtheorem{conjecture}[conjecture]{Conjecture}
\aliascntresetthe{conjecture}
\newaliascnt{problem}{theorem}

\aliascntresetthe{problem}
\newaliascnt{claim}{theorem}

\aliascntresetthe{claim}

\theoremstyle{definition}
\newaliascnt{definition}{theorem}
\newtheorem{definition}[definition]{Definition}
\aliascntresetthe{definition}
\newaliascnt{assumption}{theorem}

\aliascntresetthe{assumption}
\newaliascnt{assumptions}{theorem}

\aliascntresetthe{assumptions}
\newaliascnt{axiom}{theorem}

\aliascntresetthe{axiom}
\newaliascnt{construction}{theorem}

\aliascntresetthe{construction}
\newaliascnt{example}{theorem}
\newtheorem{example}[example]{Example}
\aliascntresetthe{example}
\newaliascnt{notation}{theorem}

\aliascntresetthe{notation}
\newaliascnt{notations}{theorem}

\aliascntresetthe{notations}
\newaliascnt{property}{theorem}

\aliascntresetthe{property}
\newaliascnt{question}{theorem}

\aliascntresetthe{question}
\newaliascnt{remark}{theorem}
\newtheorem{remark}[remark]{Remark}
\aliascntresetthe{remark}
\newaliascnt{warning}{theorem}

\aliascntresetthe{warning}
\newaliascnt{scholium}{theorem}

\aliascntresetthe{scholium}
\newaliascnt{observation}{theorem}

\aliascntresetthe{observation}
\newaliascnt{convention}{theorem}

\aliascntresetthe{convention}

\theoremstyle{plain}
\newtheorem*{claim*}{Claim}

\newtheorem{maintheorem}{Theorem}

\newaliascnt{mainproposition}{maintheorem}

\aliascntresetthe{mainproposition}
\newaliascnt{maincorollary}{maintheorem}
\newtheorem{maincorollary}[maincorollary]{Corollary}
\aliascntresetthe{maincorollary}
\newaliascnt{mainconjecture}{maintheorem}

\aliascntresetthe{mainconjecture}

\theoremstyle{definition}
\newaliascnt{mainquestion}{maintheorem}

\aliascntresetthe{mainquestion}

\numberwithin{equation}{section}
\usepackage[
  backend=biber,
  style=alphabetic,
  sorting=nyt,
  giveninits=false,
  maxbibnames=99
]{biblatex}
\usepackage{xurl} 
\AtBeginBibliography{\emergencystretch=3em}

\newif\ifhasreferences
\AtDataInput{\global\hasreferencestrue}
\DeclareFieldFormat*{title}{\mkbibemph{#1}}

\definecolor{citationgreen}{rgb}{0,0.4,0}
\usepackage[bookmarks=true]{hyperref}
\hypersetup{
  bookmarksnumbered=true,
  colorlinks=true,
  linkcolor=blue,
  citecolor=blue,
  urlcolor=blue,
  filecolor=blue,
  linktoc=page,
  pdfstartview=Fit
}

\title{Koszul duality for boundary categories}
\author{Yuan Gao}
\address{School of Mathematics, Nanjing University, Nanjing, Jiangsu, 210093, China}
\email{yuangao@nju.edu.cn} \urladdr{}
\author{Zhen Gao}
\address{Institut f\"ur Mathematik, Universit\"at Augsburg, Germany}
\email{zhen.gao@math.uni-augsburg.de} \urladdr{}
\date{\today}

\begin{document}

\begin{abstract}
We study how Koszul duality relates two kinds of boundary constructions associated with a smooth category and a proper category.
For a smooth category with a proper Koszul-dualizing subcategory, we construct a fully faithful functor from its categorical neighborhood of infinity to the singularity category associated to the proper category.
We prove that this functor is an equivalence precisely when the smooth category is reflexive, meaning that it can be reconstructed from its finite-dimensional module theory.
Under a Gorenstein condition, satisfied in the Calabi--Yau setting, the functor is induced by restriction of modules.
We give algebraic criteria that ensure equivalence and obtain an identification of Hochschild cohomology.
\end{abstract}
\maketitle
\tableofcontents

\section{Introduction}\label{sec:introduction}

\subsection{Background and motivation}

Singularity categories measure the failure of coherent complexes to be perfect and thereby make singularities accessible to homological algebra \cite{Orlov2004}.
Efimov's categorical formal punctured neighborhood of infinity provides a way to study behavior at infinity directly from a category, without choosing a compactification \cite{Efimov2017}.
A relation between these constructions allows boundary phenomena in smooth geometry to be studied through the singularities of a proper algebraic model.
Koszul duality suggests such a relation because, under suitable hypotheses, it exchanges perfect modules with modules having finite-dimensional cohomology.
To turn this principle into a comparison of boundary categories, one must identify the subcategories being divided out and determine when the comparison recovers every object of the target.
These are the two questions addressed in this paper.

The motivation is also particularly concrete in symplectic geometry.
Ganatra--Gao--Venkatesh identify the Rabinowitz Fukaya category of a nondegenerate Liouville manifold with the categorical formal punctured neighborhood of infinity of its wrapped Fukaya category \cite[Theorem 1.1]{Ganatra-Gao-Venkatesh2022}.
When compact Lagrangians supply an appropriate Koszul dual, a comparison with a singularity category offers an algebraic description of this Rabinowitz category.
The essential-surjectivity question asks whether every object of the singularity category comes from the categorical boundary, so that this algebraic description is an equivalence.

Several earlier results establish the foundations of this comparison.
Efimov relates categorical infinity to singularity categories arising from localizations of smooth proper dg categories \cite[Theorem 3.2]{Efimov2017}.
For a proper scheme over a perfect field, he identifies the categorical formal punctured neighborhood of its dg category of bounded coherent complexes with the opposite of its singularity category \cite[Theorem 9.2]{Efimov2017}.
His Proposition 9.3 gives a quotient description of infinity under categorical Koszul-duality hypotheses.
Ganatra--Gao--Venkatesh formulate these hypotheses in terms of a proper Koszul-dualizing subcategory and identify the quotient by right-proper objects with the categorical formal punctured neighborhood \cite[Definition 2.12 and Proposition 2.13]{Ganatra-Gao-Venkatesh2022}.
We use this quotient description to compare the boundary with the singularity category of the Koszul dual.

At the level of modules, the exchange of finiteness conditions is a familiar feature of Koszul duality.
For Adams connected $A_\infty$ algebras with proper Koszul dual, Lu--Palmieri--Wu--Zhang prove an equivalence between perfect modules on one side and modules with finite-dimensional total cohomology on the other \cite[Theorem B and Corollary 7.2(b)]{Lu-Palmieri-Wu-Zhang2008}.
The question of reconstructing a category from its pseudo-perfect modules is expressed more generally by \emph{reflexivity}, introduced by Kuznetsov--Shinder \cite[Definition 3.11]{Kuznetsov-Shinder2025}.
Goodbody characterizes reflexivity in the Morita homotopy category and studies its consequences for Hochschild cohomology \cite[Theorems 4.4 and 5.1]{Goodbody2026}.
Booth--Goodbody--Opper develop its relationship with Koszul duality, generation, and completion, and establish criteria for connective and coconnective dg algebras \cite[Theorem 2.3.8 and Section 3]{Booth-Goodbody-Opper2025}.
These results provide the reconstruction property and the algebraic criteria used below.

Calabi--Yau structures give a further reason to compare the two boundaries.
Brav--Dyckerhoff transfer a smooth Calabi--Yau structure to the subcategory of proper objects \cite[Theorem 3.1]{Brav-Dyckerhoff2019}, while Gao's categorical residue produces a duality of dimension one less on the categorical formal punctured neighborhood of infinity \cite[Theorem 1.7]{Gao2025}.
Rivera--Takeda--Wang study the associated products on Hochschild complexes with boundary coefficients using pre-Calabi--Yau structures \cite{Rivera-Takeda-Wang2024}.
In his Oberwolfach report, Takeda conjectures compatibility of the products and duality structures appearing in Rabinowitz Floer theory, categorical infinity, and singular Hochschild cohomology under the relevant identifications \cite[p. 3165]{Takeda2023}, which leads to the general form of the following conjecture:

\begin{conjecture}\label{conjecture:Takeda}
Let $\mathscr{B}$ be a proper Koszul dualizing subcategory in a smooth $A_\infty$-category $\mathscr{A}$. Then there is an isomorphism of rings:
\[
HH_{\mathsf{sg}}^*(\mathscr{B})\cong HH_*(\widehat{\mathscr{A}}_{\infty}) 
\]
\end{conjecture}

The main purpose of this paper is to address the categorical identification needed for this perspective and give consequences for Hochschild invariants.

\subsection{Main results}

Fix a field $\Bbbk$.
All starting categories are small and $\Bbbk$-linear, and we use dg module categories, passing to dg models for strictly unital $A_\infty$ categories when forming quotients.
All tensor products and module Hom functors are computed in these models.
For a category $\mathscr C$, write $\mathsf{Mod}(\mathscr C)$ for its dg category of right modules, $\mathsf{Perf}(\mathscr C)$ for the full subcategory of perfect modules, and $\mathsf{PsPerf}(\mathscr C)$ for the full subcategory of pseudo-perfect modules, whose values are perfect $\Bbbk$-complexes.
Smoothness gives $\mathsf{PsPerf}(\mathscr C)\subset\mathsf{Perf}(\mathscr C)$, while properness gives the reverse inclusion.
For proper $\mathscr B$, the dg singularity category is
\[
\mathscr D_{\mathsf{sg}}(\mathscr B)
=\mathsf{PsPerf}(\mathscr B)/\mathsf{Perf}(\mathscr B).
\]

Let $j:\mathscr B\hookrightarrow\mathscr A$ be a Koszul-dualizing full inclusion in the sense of Definition \ref{def:Koszul-duality}, with $\mathscr A$ smooth and $\mathscr B$ proper.
In particular, the graph bimodule $\Gamma_j(a,b)=\mathrm{Hom}_{\mathscr A}(a,jb)$ is perfect over $\Bbbk$ at every pair of objects.
Writing $(-)^\vee=\mathrm{Hom}_{\Bbbk}(-,\Bbbk)$ for linear duality, define
\begin{equation}\label{eq:compare-functor}
\Phi=-\otimes_{\mathscr A}\Gamma_j^\vee:
\mathsf{Perf}(\mathscr A)\longrightarrow\mathsf{PsPerf}(\mathscr B).
\end{equation}
For a perfect module $M$, its value at $b$ is $\Phi(M)(b)\simeq\mathrm{Hom}_{\mathscr A}(M,\mathbf Y_{\mathscr A}^{\mathrm r}(jb))^\vee$, where $\mathbf Y_{\mathscr A}^{\mathrm r}(jb)=\mathrm{Hom}_{\mathscr A}(-,jb)$ is the right Yoneda module.
Equivalently,
\begin{equation}\label{eq:compare-functor-Hom-form}
\Phi\simeq(-)^\vee\circ\bigl((j^{\mathrm{op}})^*\bigr)^{\mathrm{op}}
\circ\bigl(\mathrm{Hom}_{\mathscr A}(-,\mathscr A)\bigr)^{\mathrm{op}},
\end{equation}
where the rightmost functor is perfect-module duality and $(j^{\mathrm{op}})^*$ is restriction of left modules.

We will need the notion of \textbf{reflexivity}.
The category $\mathscr A$ is \emph{reflexive} if the canonical evaluation functor
\[
\begin{aligned}
\operatorname{ev}_{\mathscr A}:
\mathsf{Perf}(\mathscr A)^{\mathrm{op}}
&\longrightarrow\mathsf{PsPerf}\bigl(\mathsf{PsPerf}(\mathscr A)^{\mathrm{op}}\bigr),\\
M&\longmapsto\bigl[N\longmapsto\mathrm{Hom}_{\mathscr A}(M,N)\bigr]
\end{aligned}
\]
is a quasi-equivalence; see Definition \ref{def:reflexivity} for further details.
Write $\widehat{\mathsf{Perf}(\mathscr A)}_\infty$ for the categorical formal punctured neighborhood of infinity of $\mathsf{Perf}(\mathscr A)$, using the essential-image construction of Definition \ref{def:cat-form-punc-nbhd}.
Under the present hypotheses, it is quasi-equivalent to $\mathsf{Perf}(\mathscr A)/\mathsf{PsPerf}(\mathscr A)$ by Proposition \ref{prop:boundary-quotient}.

\begin{maintheorem}\label{thm:main}
Let $\mathscr A$ be smooth and let $j:\mathscr B\hookrightarrow\mathscr A$ be a Koszul-dualizing full inclusion, with $\mathscr B$ proper.
The functor $\Phi$ in \eqref{eq:compare-functor} is quasi-fully faithful.
Its composite with linear duality induces a quasi-fully faithful functor
\begin{equation}\label{eq:main-comparison}
\widehat{\mathsf{Perf}(\mathscr A)}_\infty
\longrightarrow\mathscr D_{\mathsf{sg}}(\mathscr B^{\mathrm{op}})^{\mathrm{op}}.
\end{equation}
The following conditions are equivalent:
\begin{enumerate}
    \item[(i)] The functor \eqref{eq:main-comparison} is a quasi-equivalence.
    \item[(ii)] The functor $\Phi$ is a quasi-equivalence.
    \item[(iii)] The category $\mathscr A$ is reflexive.
\end{enumerate}
If $\mathscr B$ also satisfies the Gorenstein condition \eqref{eq:Gorenstein-cond}, its Nakayama functor $\nu_{\mathscr B}=-\otimes_{\mathscr B}\mathscr B^\vee$, where $\mathscr B^\vee$ is the dual diagonal bimodule, is an autoequivalence of $\mathsf{Mod}(\mathscr B)$ preserving $\mathsf{PsPerf}(\mathscr B)$, and
\[
\Phi\simeq\nu_{\mathscr B}\circ j^*
\qquad\text{on }\mathsf{Perf}(\mathscr A).
\]
Restriction $j^*$ is then quasi-fully faithful on $\mathsf{Perf}(\mathscr A)$ and induces a quasi-fully faithful functor
\begin{equation}\label{eq:ordinary-comparison}
\widehat{\mathsf{Perf}(\mathscr A)}_\infty
\longrightarrow\mathscr D_{\mathsf{sg}}(\mathscr B).
\end{equation}
This functor is a quasi-equivalence precisely when \textup{(i)--(iii)} hold.
These conditions are also equivalent to coinduction, the right adjoint of $j^*$, restricting to a quasi-equivalence
\[
j_*:\mathsf{PsPerf}(\mathscr B)\xrightarrow{\;\sim\;}\mathsf{Perf}(\mathscr A),
\]
with quasi-inverse $j^*$.
\end{maintheorem}

The proof separates reconstruction from passage to the quotient.
Efimov's generation result identifies $\mathsf{Perf}(\mathscr B)$ with $\mathsf{PsPerf}(\mathscr A)$; under this identification, the factorization \eqref{eq:compare-functor-Hom-form} expresses $\Phi$ through the evaluation functor defining reflexivity.
For the quotients, $\Phi$ sends $\mathsf{PsPerf}(\mathscr A)$ to the Nakayama image $\nu_{\mathscr B}(\mathsf{Perf}(\mathscr B))$, which linear duality identifies with $\mathsf{Perf}(\mathscr B^{\mathrm{op}})^{\mathrm{op}}$.
This determines the target in \eqref{eq:main-comparison}.
The essential image of $\Phi$ is thick, so quotienting by the corresponding full subcategories preserves and detects essential surjectivity, by Lemma \ref{lem:quotient-density}.
The Gorenstein condition identifies the Nakayama image with $\mathsf{Perf}(\mathscr B)$ and gives the restriction functor in \eqref{eq:ordinary-comparison}.

A weak proper $n$-Calabi--Yau structure on $\mathscr B$ identifies $\mathscr B^\vee$ with $\mathscr B_\Delta[n]$.
It therefore implies the Gorenstein condition and gives $\nu_{\mathscr B}\simeq[n]$, so that $\Phi\simeq[n]\circ j^*$.
Combining Theorem \ref{thm:main} with the reflexivity criteria of Booth--Goodbody--Opper yields the following consequence.
The degreewise finiteness required by those criteria follows from the Koszul-duality hypotheses: the tensor--Hom map in Definition \ref{def:Koszul-duality}(iii) is a biduality map, as in \eqref{eq:algebraic-biduality}, and condition (ii) transfers this finiteness to morphisms between perfect $\mathscr A$-modules.

\begin{maincorollary}\label{cor:main}
Let $\mathscr A$ be smooth and let $j:\mathscr B\hookrightarrow\mathscr A$ be a Koszul-dualizing full inclusion.
Suppose that $\mathscr B$ carries a weak proper $n$-Calabi--Yau structure and that $\mathsf{Perf}(\mathscr A)$ has a split-generator $G$ whose dg endomorphism algebra $A=\mathrm{Hom}_{\mathscr A}(G,G)$ satisfies either
\begin{equation}\label{eq:connective}
H^{>0}(A)=0,
\end{equation}
or
\begin{equation}\label{eq:coconnective}
H^{<0}(A)=0
\qquad\text{and}\qquad
H^0(A)\text{ is semisimple}.
\end{equation}
Then restriction induces a quasi-equivalence
\[
\widehat{\mathsf{Perf}(\mathscr A)}_\infty
\xrightarrow{\;\sim\;}\mathscr D_{\mathsf{sg}}(\mathscr B).
\]
Moreover, $j^*$ and $j_*$ restrict to mutually inverse quasi-equivalences between $\mathsf{Perf}(\mathscr A)$ and $\mathsf{PsPerf}(\mathscr B)$.
\end{maincorollary}

If $\mathscr A$ itself carries a weak smooth $n$-Calabi--Yau structure, Theorem \ref{thm:smooth-proper-CY} gives the weak proper structure on $\mathscr B$.
When $\mathscr A$ is reflexive, Theorem \ref{thm:koszul-CY-comparison} transports Gao's categorical residue through the restriction quasi-equivalence and equips $\mathscr D_{\mathsf{sg}}(\mathscr B)$ with a weak generalized proper $(n-1)$-Calabi--Yau structure in the sense of Remark \ref{rem:CY-finiteness}.
Corollary \ref{cor:HH-boundary} gives the corresponding isomorphisms of Hochschild cohomology algebras, preserving cup products.

\subsection*{Structure of the paper}

In Section \ref{sec:basic}, we recall module categories, the Nakayama functor, and the properties of dg quotients used below.
In Section \ref{sec:calabi-yau}, we review weak smooth and weak proper Calabi--Yau structures and their relation to the Gorenstein condition.
In Section \ref{sec:koszul-reflexivity}, we prove that the comparison functor $\Phi$ is quasi-fully faithful and characterize its essential surjectivity by reflexivity.
In Section \ref{sec:infinity-singularity}, we recall the two boundary constructions and prove Theorem \ref{thm:main} and its Calabi--Yau refinement, Theorem \ref{thm:koszul-CY-comparison}.
In Section \ref{sec:algebraic-criteria}, we give algebraic criteria for the comparison to be a quasi-equivalence, apply them to Koszul duals of proper algebras, and discuss geometric examples.
In Section \ref{sec:hochschild}, we deduce comparisons of Hochschild cohomology algebras and relate them to singular Hochschild cohomology under Keller's hypotheses.

\subsection*{Acknowledgments}

We thank Alex Takeda for sharing his inspiring perspective, and Zhengfang Wang for helpful discussion, especially for suggesting the \textit{reflexivity} property for consideration.

\section{Categorical preliminaries}\label{sec:basic}

We recall the finiteness conditions, duality functors, and properties of dg quotients.
All starting categories are small (except for module categories) and $\Bbbk$-linear, complexes are cohomologically graded, and modules are right modules unless otherwise specified.
Tensor products and Hom functors are derived.
We use the standard module and Morita theory of dg and strictly unital $A_\infty$ categories; see \cite{Keller2006a,Keller2006b,Toen2007,Seidel2008}.

\subsection{Basic definitions}
\begin{definition}\label{def:perf-prop}
	Let $\mathscr{C}$ be a small dg or strictly unital $A_\infty$ category.
	\begin{enumerate}
		\item $\mathsf{Mod}(\mathscr{C})$ is the dg category of right $\mathscr{C}$-modules, defined as the dg category of dg or $A_\infty$ functor $\mathsf{Fun}(\mathscr{C}^{\mathrm{op}},\mathsf{Mod}(\Bbbk))$.
		
		\item $\mathsf{Perf}(\mathscr{C})$ is the split-closed pretriangulated envelope of
		the right Yoneda modules $\mathbf{Y}_{\mathscr{C}}^{\mathrm{r}}(X):=\mathrm{Hom}_{\mathscr{C}}(-,X)$.  Equivalently, it is the full subcategory of compact objects of $\mathsf{Mod}(\mathscr{C})$, i.e. $\mathsf{Perf}(\mathscr{C})=\mathsf{Mod}(\mathscr{C})^{\omega}$.
		
		\item $\mathsf{PsPerf}(\mathscr{C})$ is the full subcategory of modules $M$ for which $M(X)$ is a perfect $\Bbbk$-complex for every $X\in\mathscr{C}$. Such modules are called or \textbf{pseudo-perfect} or \textbf{proper}.
		
		\item An object $X\in\mathscr{C}$ is \textbf{right proper} if its right Yoneda module $\mathbf{Y}_{\mathscr{C}}^{\mathrm{r}}(X)$ is proper.  The full subcategory of right-proper objects is denoted $\mathscr{C}^{\mathrm{prop}}$.
	\end{enumerate}
\end{definition}

There are following quasi-equivalences of module categories:
\begin{eqnarray}\label{eq:op-equivalence}
\mathrm{Hom}_{\mathscr{C}}(-,\mathscr{C}) : \mathsf{Perf}(\mathscr{C})^{\mathrm{op}} &\xrightarrow{\;\sim\;} \mathsf{Perf}(\mathscr{C}^{\mathrm{op}}) , \label{eq:perf-op-equivalence}\\ 
(-)^\vee : \mathsf{PsPerf}(\mathscr{C})^{\mathrm{op}} &\xrightarrow{\;\sim\;} \mathsf{PsPerf}(\mathscr{C}^{\mathrm{op}}) .  \label{eq:psperf-op-equivalence}
\end{eqnarray}

\begin{remark}
For a strictly unital $A_\infty$ category, the Yoneda embedding gives an $A_\infty$ quasi-equivalence from $\mathscr C$ to the full dg subcategory of $\mathsf{Mod}(\mathscr C)$ on its representable modules \cite[Sections 2f--2g]{Seidel2008}.
We use these dg models when forming tensor products of dg categories and Drinfeld quotients.
\end{remark}

Set $\mathscr C^{\mathrm e}=\mathscr C^{\mathrm{op}}\otimes_{\Bbbk}\mathscr C$, and use the diagonal convention
\[
\mathscr C_\Delta(X,Y)=\mathrm{Hom}_{\mathscr C}(Y,X).
\]
Thus the \textbf{diagonal bimodule} $\mathscr C_\Delta$ is a right $\mathscr C^{\mathrm e}$-module, covariant in $X$ and contravariant in $Y$.

\begin{definition}\label{def:smooth-proper}
The category $\mathscr C$ is \textbf{smooth} if $\mathscr C_\Delta\in\mathsf{Perf}(\mathscr C^{\mathrm e})$.
It is \textbf{proper} if $\mathrm{Hom}_{\mathscr C}(X,Y)\in\mathsf{Perf}(\Bbbk)$ for all $X,Y\in\mathscr C$, or equivalently if $\mathscr C_\Delta\in\mathsf{PsPerf}(\mathscr C^{\mathrm e})$.
Over the field $\Bbbk$, the latter condition means that each morphism complex has finite-dimensional total cohomology.
\end{definition}

The following elementary facts relate these finiteness conditions and will also be used to check preservation of pseudo-perfect modules.

\begin{lemma}\label{lem:module-finiteness}
A module $M$ is pseudo-perfect if and only if $\mathrm{Hom}_{\mathscr C}(P,M)\in\mathsf{Perf}(\Bbbk)$ for every $P\in\mathsf{Perf}(\mathscr C)$.
Moreover,
\[
\begin{aligned}
\mathscr C\text{ proper}&\quad\Longrightarrow\quad\mathsf{Perf}(\mathscr C)\subset\mathsf{PsPerf}(\mathscr C),\\
\mathscr C\text{ smooth}&\quad\Longrightarrow\quad\mathsf{PsPerf}(\mathscr C)\subset\mathsf{Perf}(\mathscr C).
\end{aligned}
\]
\end{lemma}
\begin{proof}
The first assertion follows from $\mathrm{Hom}_{\mathscr C}(\mathbf Y_{\mathscr C}^{\mathrm r}(X),M)\simeq M(X)$ and closure of $\mathsf{Perf}(\Bbbk)$ under finite cones and direct summands.
If $\mathscr C$ is proper, its right Yoneda modules are pseudo-perfect, and the same closure argument gives the first inclusion.
If $\mathscr C$ is smooth, tensor a finite resolution of $\mathscr C_\Delta$ from representable bimodules with a pseudo-perfect module $M$.
Since
\[
M\otimes_{\mathscr C} 
\bigl(\mathbf Y_{\mathscr C}^{\mathrm l}(X)\otimes_{\Bbbk}\mathbf Y_{\mathscr C}^{\mathrm r}(Y)\bigr)
\simeq M(X)\otimes_{\Bbbk}\mathbf Y_{\mathscr C}^{\mathrm r}(Y)
\]
is perfect, so is $M\simeq M\otimes_{\mathscr C} \mathscr C_\Delta$.
\end{proof}

\subsection{The Nakayama functor}

The \textbf{dual diagonal bimodule} is defined to be $\mathscr C^\vee(X,Y)=\mathrm{Hom}_{\mathscr C}(X,Y)^\vee$.
The Nakayama functor is
\begin{equation}\label{eq:Nakayama-functor}
\nu_{\mathscr C}=-\otimes_{\mathscr C} \mathscr C^\vee:
\mathsf{Mod}(\mathscr C)\longrightarrow\mathsf{Mod}(\mathscr C).
\end{equation}
For $P\in\mathsf{Perf}(\mathscr C)$, tensor--Hom duality gives
\[
\nu_{\mathscr C}(P)\simeq\mathrm{Hom}_{\mathscr C}(P,\mathscr C)^\vee.
\]
Its right adjoint is
\begin{equation}\label{eq:Nakayama-functor-right-adjoint}
\nu_{\mathscr{C}}^{\mathsf{R}}=\mathrm{Hom}_{\mathscr{C}}(\mathscr{C}^\vee,-) : \mathsf{Mod}(\mathscr{C}) \longrightarrow \mathsf{Mod}(\mathscr{C}).
\end{equation}

Suppose now that $\mathscr C$ is proper, and put
\begin{equation}\label{eq:Nakayama-image}
\mathscr I_{\mathscr C}
=\mathrm{thick}_{\mathsf{Mod}(\mathscr C)}
\{\,\mathbf Y_{\mathscr C}^{\mathrm l}(X)^\vee:X\in\mathscr C\,\}
\subset\mathsf{PsPerf}(\mathscr C).
\end{equation}
This is the image of the perfect modules under the Nakayama functor, which we will show in Lemma \ref{lem:Nakayama-Serre}.
We also recall the following definition.

\begin{definition}\label{def:Gorenstein}
A proper category $\mathscr C$ satisfies the \textbf{Gorenstein condition} if
\begin{equation}\label{eq:Gorenstein-cond}
\mathscr I_{\mathscr C}=\mathsf{Perf}(\mathscr C)
\end{equation}
as full subcategories of $\mathsf{Mod}(\mathscr C)$.
For proper dg algebras, this is the condition of \cite[Definition 5.2]{Goodbody2024}.
\end{definition}

\begin{remark}
For an ordinary finite-dimensional algebra $B$, one has $\mathscr I_B=\mathrm{thick}(B^\vee)=\mathsf K^{\mathrm b}(\mathrm{inj}\text{-}B)$, where $\mathrm{inj}\text{-}B$ denotes the finitely generated injective right modules.
Gorenstein ondition \eqref{eq:Gorenstein-cond} is equivalent to finite injective dimension of both regular modules ${}_BB$ and $B_B$; see \cite[Corollary 3.9]{Chen2011} and \cite[Remark 5.3]{Goodbody2024}.
\end{remark}

\begin{lemma}\label{lem:Nakayama-Serre}
If $\mathscr C$ is proper, the Nakayama functor restricts to a quasi-equivalence
\[
\nu_{\mathscr C}:\mathsf{Perf}(\mathscr C)\xrightarrow{\;\sim\;}\mathscr I_{\mathscr C}.
\]
For $P\in\mathsf{Perf}(\mathscr C)$ and $M\in\mathsf{Mod}(\mathscr C)$, there is a natural quasi-isomorphism
\begin{equation}\label{eq:Nakayama-duality}
\mathrm{Hom}_{\mathscr C}(P,M)^\vee
\simeq\mathrm{Hom}_{\mathscr C}(M,\nu_{\mathscr C}P).
\end{equation}
If $\mathscr C$ satisfies the Gorenstein condition, then $\nu_{\mathscr C}$ is a Serre functor on $\mathsf{Perf}(\mathscr C)$ and an autoequivalence of $\mathsf{Mod}(\mathscr C)$ preserving $\mathsf{PsPerf}(\mathscr C)$.
Its inverse is $\nu_{\mathscr C}^{\mathsf R}$.
\end{lemma}
\begin{proof}
On perfect modules, $\nu_{\mathscr C}$ is the composite of the dualities in \eqref{eq:op-equivalence}, using the inclusion $\mathsf{Perf}(\mathscr C^{\mathrm{op}})\subset\mathsf{PsPerf}(\mathscr C^{\mathrm{op}})$ since $\mathscr C$ is proper.
It is therefore quasi-fully faithful.
Since it sends $\mathbf Y_{\mathscr C}^{\mathrm r}(X)$ to $\mathbf Y_{\mathscr C}^{\mathrm l}(X)^\vee$ and its source is idempotent complete, its essential image is precisely $\mathscr I_{\mathscr C}$.
For representable $P$, \eqref{eq:Nakayama-duality} is the tensor--Hom adjunction
\[
\mathrm{Hom}_{\mathscr C}
\bigl(M,\mathbf Y_{\mathscr C}^{\mathrm l}(X)^\vee\bigr)
\simeq\bigl(M\otimes_{\mathscr C} \mathbf Y_{\mathscr C}^{\mathrm l}(X)\bigr)^\vee
\simeq M(X)^\vee;
\]
the general case follows by finite cones and direct summands.

Under \eqref{eq:Gorenstein-cond}, this gives a Serre autoequivalence of $\mathsf{Perf}(\mathscr C)$.
The bar tensor functor $\nu_{\mathscr C}$ preserves direct sums and cones and induces a quasi-equivalence on the compact objects, so its colimit-preserving extension gives an autoequivalence of $\mathsf{Mod}(\mathscr C)$ \cite[Lemma 3.10]{Keller2006b}.
By Lemma \ref{lem:module-finiteness}, pseudo-perfectness is characterized by perfectness over $\Bbbk$ of the Hom complexes from all perfect modules.
This condition is preserved by $\nu_{\mathscr C}$ and its inverse, proving the last assertion.
\end{proof}

\subsection{Drinfeld quotient}

Let $\mathscr N\subset\mathscr C$ be a full dg subcategory.
The \textbf{dg quotient} is a dg category $\mathscr C/\mathscr N$ equipped with a functor
\[
q:\mathscr C\longrightarrow\mathscr C/\mathscr N
\]
that annihilates $\mathscr N$ and is universal with this property in the homotopy category of dg categories \cite[Theorem 1.6.2]{Drinfeld2004}, \cite[Theorem 4.8]{Keller2006b}.
More precisely, precomposition with $q$ identifies morphisms $\mathscr C/\mathscr N\to\mathscr D$ in $\mathsf{Ho}(\mathsf{dgCat}_{\Bbbk})$ with morphisms $\mathscr C\to\mathscr D$ that send every object of $\mathscr N$ to a contractible object.
Over the field $\Bbbk$, a concrete model has the same objects as $\mathscr C$ and freely adjoins, for each $N\in\mathscr N$, an endomorphism $\varepsilon_N$ of degree $-1$ with $d\varepsilon_N=\mathrm{id}_N$ \cite[Section 3.1]{Drinfeld2004}.
Its morphism complexes admit the description
\begin{equation}\label{eq:quotient-Hom}
\mathrm{Hom}_{\mathscr C/\mathscr N}(X,Y)
\simeq\mathrm{cone}\left(
\mathbf Y_{\mathscr C}^{\mathrm r}(Y)|_{\mathscr N}
\otimes_{\mathscr N} 
\mathbf Y_{\mathscr C}^{\mathrm l}(X)|_{\mathscr N}
\longrightarrow\mathrm{Hom}_{\mathscr C}(X,Y)
\right),
\end{equation}
where the arrow is induced by composition \cite[Section 1.3.1]{Drinfeld2004}.

If $\mathscr C$ and $\mathscr N$ are pretriangulated, then $\mathscr C/\mathscr N$ is pretriangulated and there is an equivalence of triangulated categories
\[
H^0(\mathscr C/\mathscr N)\simeq H^0(\mathscr C)/H^0(\mathscr N)
\]
with the Verdier quotient \cite[Theorem 3.4]{Drinfeld2004}.
In particular,
\[
q(X)\simeq0
\quad\Longleftrightarrow\quad
X\in\mathrm{thick}_{H^0(\mathscr C)}H^0(\mathscr N)
\]
by \cite[Remark A.3(ii)]{Drinfeld2004} or \cite[Proposition 4.6.2(3)]{Krause2010}.
Consequently, replacing $\mathscr N$ by its thick closure inside $\mathscr C$ gives a quasi-equivalent dg quotient.
Quasi-equivalences of dg pairs also induce quasi-equivalences of their quotients \cite[Remark 3.6(ii)]{Drinfeld2004}.
For $A_\infty$ categories, one may use dg models or the quotient constructions of \cite{Lyubashenko-Ovsienko2006,Lyubashenko-Manzyuk2008}.

The following lemma shows that quotienting by a common subcategory preserves full faithfulness and, under a thickness hypothesis, detects essential surjectivity; see \cite[Lemma A.5]{Drinfeld2004},\cite[Proposition 4.6.2, Lemma 4.7.1]{Krause2010} for the arguments for triangulated categories.

\begin{lemma}\label{lem:quotient-density}
Let $\mathscr N\subset\mathscr C\subset\mathscr D$ be full pretriangulated dg subcategories.
The induced functor
\[
\bar j:\mathscr C/\mathscr N\longrightarrow\mathscr D/\mathscr N
\]
is quasi-fully faithful.
If $H^0(\mathscr C)$ is thick in $H^0(\mathscr D)$, then $\bar j$ is essentially surjective if and only if the inclusion $j:\mathscr C\hookrightarrow\mathscr D$ is essentially surjective.
\end{lemma}
\begin{proof}
For $X,Y\in\mathscr C$, the cone descriptions \eqref{eq:quotient-Hom} in $\mathscr C/\mathscr N$ and $\mathscr D/\mathscr N$ agree because both inclusions are full.
Thus $\bar j$ is quasi-fully faithful.

Suppose that $H^0(\mathscr C)$ is thick in $H^0(\mathscr D)$.
The quotient $\mathscr D\to\mathscr D/\mathscr C$ factors through $\mathscr D/\mathscr N$.
By the kernel description above, an object represented by $X\in\mathscr D$ vanishes in $H^0(\mathscr D/\mathscr C)$ precisely when $X$ is isomorphic in $H^0(\mathscr D)$ to an object of $\mathscr C$.
Hence
\begin{equation}\label{eq:kernel}
\operatorname{ess\,im}H^0(\bar j)
=\ker\bigl(H^0(\mathscr D/\mathscr N)\longrightarrow H^0(\mathscr D/\mathscr C)\bigr).
\end{equation}
If $\bar j$ is essentially surjective, every object of $\mathscr D$ therefore belongs to $H^0(\mathscr C)$ up to isomorphism.
Conversely, essential surjectivity of $j$ passes to the quotients, which have the same objects as their respective source categories.
\end{proof}

\begin{remark}
Under the thickness hypothesis of Lemma \ref{lem:quotient-density}, the kernel as in \eqref{eq:kernel} shows that the essential image of $H^0(\bar j)$ is thick.
A quasi-fully faithful functor between pretriangulated dg categories is a Morita equivalence precisely when its image generates the target under shifts, cones, and direct summands \cite[Lemma 3.10]{Keller2006b}.
Thus $\bar j$ is a Morita equivalence if and only if it is a quasi-equivalence.
\end{remark}

\section{Calabi--Yau structures}\label{sec:calabi-yau}

We recall weak smooth and weak proper Calabi--Yau structures, which will be used in Theorem \ref{thm:koszul-CY-comparison} and Corollary \ref{cor:main}, respectively.
Lemma \ref{lem:pCY-Gorenstein} shows that a weak proper Calabi--Yau structure implies the Gorenstein condition \eqref{eq:Gorenstein-cond} required in Theorem \ref{thm:main} for the functor \eqref{eq:ordinary-comparison}.

Write $C_*(\mathscr C)=C_*(\mathscr C,\mathscr C_\Delta)$ for the Hochschild chain complex with diagonal coefficients, computed by the cyclic bar construction, and set
\[
HH_m(\mathscr C)=H^{-m}(C_*(\mathscr C)),\qquad m\in\mathbb Z.
\]
We regard $C_*(\mathscr C)$ as a cochain complex; the chain maps defining Calabi--Yau structures have degree zero.
Thus a class in $HH_n(\mathscr C)$ is represented by a chain map $\Bbbk[n]\to C_*(\mathscr C)$.
For the cyclic bar construction, see \cite[Section 2.4]{Gao2025}.

The definition of a weak smooth Calabi--Yau structure uses the \textbf{inverse dualizing bimodule} $\mathscr C^!$, the $\mathscr C^{\mathrm e}$-linear dual of $\mathscr C_\Delta$, regarded as a right $\mathscr C^{\mathrm e}$-module by exchanging the factors of $\mathscr C^{\mathrm e}=\mathscr C^{\mathrm{op}}\otimes_{\Bbbk}\mathscr C$.
Explicitly,
\[
\mathscr C^!(X,Y)=\mathrm{Hom}_{\mathscr C^{\mathrm e}}
\bigl(\mathscr C_\Delta,
\mathbf Y_{\mathscr C}^{\mathrm l}(Y)\otimes_{\Bbbk}\mathbf Y_{\mathscr C}^{\mathrm r}(X)\bigr).
\]
The tensor product in the second argument is the bimodule whose value at $(U,V)$ is $\mathrm{Hom}_{\mathscr C}(Y,U)\otimes_{\Bbbk}\mathrm{Hom}_{\mathscr C}(V,X)$.
If $\mathscr C$ is smooth, perfect duality for $\mathscr C_\Delta$ gives a natural quasi-isomorphism
\begin{equation}\label{eq:CY-smooth-Hochschild}
C_*(\mathscr C)\simeq\mathrm{Hom}_{\mathscr C^{\mathrm e}}(\mathscr C^!,\mathscr C_\Delta);
\end{equation}
see \cite[Section 3.2]{Brav-Dyckerhoff2019} and \cite[Remark 6.7]{Ganatra2023}.

\begin{definition}\label{def:weak-smooth-CY}
Suppose that $\mathscr C$ is smooth.
A \textbf{weak smooth $n$-Calabi--Yau structure} is a class $[\omega_{\mathscr C}]\in HH_n(\mathscr C)$, represented by
\begin{equation}\label{eq:sCY-volume}
\omega_{\mathscr C}:\Bbbk[n]\longrightarrow C_*(\mathscr C),
\end{equation}
for which the bimodule morphism corresponding under \eqref{eq:CY-smooth-Hochschild} is a quasi-isomorphism
\begin{equation}\label{eq:sCY-map}
\mathscr C^![n]\xrightarrow{\,\sim\,}\mathscr C_\Delta.
\end{equation}
See \cite[Definition 6.5 and Remark 6.7]{Ganatra2023}.
\end{definition}

The tensor--Hom adjunction gives a natural quasi-isomorphism
\begin{equation}\label{eq:CY-proper-Hochschild}
\mathrm{Hom}_{\Bbbk}(C_*(\mathscr C),\Bbbk)
\simeq\mathrm{Hom}_{\mathscr C^{\mathrm e}}(\mathscr C_\Delta,\mathscr C^\vee),
\end{equation}
where the right-hand side is the morphism complex in $\mathsf{Mod}(\mathscr C^{\mathrm e})$; see \cite[Section 3.1]{Brav-Dyckerhoff2019} and \cite[Remark 6.2]{Ganatra2023}.

\begin{definition}\label{def:weak-proper-CY}
Suppose that $\mathscr C$ is proper.
A \textbf{weak proper $n$-Calabi--Yau structure} is a chain map, called a Hochschild trace,
\begin{equation}\label{eq:pCY-trace}
\operatorname{tr}:C_*(\mathscr C)\longrightarrow\Bbbk[-n],
\end{equation}
taken up to chain homotopy, for which the bimodule morphism corresponding under \eqref{eq:CY-proper-Hochschild} is a quasi-isomorphism
\begin{equation}\label{eq:pCY-map}
\mathscr C_\Delta\xrightarrow{\,\sim\,}\mathscr C^\vee[-n].
\end{equation}
For a trace as in \eqref{eq:pCY-trace}, this quasi-isomorphism condition is equivalent to requiring that composition followed by the trace on length-zero Hochschild chains induce perfect pairings
\begin{equation}\label{eq:CY-pairing}
H^i\mathrm{Hom}_{\mathscr C}(X,Y)\otimes_{\Bbbk}
H^{n-i}\mathrm{Hom}_{\mathscr C}(Y,X)\longrightarrow\Bbbk
\end{equation}
for every $X,Y\in\mathscr C$ and $i\in\mathbb Z$; see \cite[Definition 6.1 and Remark 6.2]{Ganatra2023}.
\end{definition}

\begin{remark}\label{rem:CY-finiteness}
In Definition \ref{def:weak-proper-CY}, one can drop the a priori condition of properness on $\mathscr{C}$. The perfect pairing, induced by the Hochschild trace \eqref{eq:pCY-trace} which determines the bimodule quasi-isomorphism \eqref{eq:pCY-map} via the derived tensor--Hom adjunction
\[
\mathrm{Hom}_{\Bbbk}(C_*(\mathscr{C}),\Bbbk[-n])\simeq \mathrm{Hom}_{\mathscr{C}^{\mathrm{e}}}(\mathscr{C}_{\Delta},\mathscr{C}^\vee[-n]),
\]
implies the finiteness of Hom spaces on homology level at each degree 
\[
H^i(\mathrm{Hom}_{\mathscr{C}}(X,Y))\xrightarrow{\;\sim\;}H^{n-i}(\mathrm{Hom}_{\mathscr{C}}(Y,X))^\vee
\]
(whereas the total dimension of homology can nevertheless be infinite), and this means $\mathscr{C}$ is a locally proper dg or $A_\infty$ category. In this spirit, we refer to this milder notion as \textbf{weak (absolute) locally proper Calabi--Yau structure} on $\mathscr{C}$; this is also named a \textbf{component-wise proper Calabi--Yau structure} by Keller--Wang.
There is an even weaker notion in \cite[Definition 2.27]{Gao2025}, obtained by requiring that the pairings \eqref{eq:CY-pairing} be nondegenerate rather than perfect.
The corresponding bimodule morphism $\mathscr{C}_{\Delta}\to\mathscr{C}^\vee[-n]$ is then injective on cohomology, with no surjectivity requirement, and the morphism cohomology groups may be infinite-dimensional in each degree.
We call such a structure a \textbf{weak generalized proper $n$-Calabi--Yau structure}.
\end{remark} 

The next lemma verifies the Gorenstein hypothesis of Theorem \ref{thm:main} for the category $\mathscr B$ in Corollary \ref{cor:main}.

\begin{lemma}\label{lem:pCY-Gorenstein}
If $\mathscr C$ carries a weak proper $n$-Calabi--Yau structure, then $\nu_{\mathscr C}\simeq[n]$ on $\mathsf{Mod}(\mathscr C)$.
In particular, $\mathscr C$ satisfies the Gorenstein condition \eqref{eq:Gorenstein-cond}, and the Serre functor on $\mathsf{Perf}(\mathscr C)$ is the shift $[n]$.
\end{lemma}
\begin{proof}
The quasi-isomorphism \eqref{eq:pCY-map} gives $\mathscr C^\vee\simeq\mathscr C_\Delta[n]$.
Substituting this identification into the definition \eqref{eq:Nakayama-functor} of the Nakayama functor yields
\[
\nu_{\mathscr C}
=-\otimes_{\mathscr C} \mathscr C^\vee
\simeq-\otimes_{\mathscr C} \mathscr C_\Delta[n]
\simeq[n].
\]
By Lemma \ref{lem:Nakayama-Serre}, the essential image of $\nu_{\mathscr C}$ on $\mathsf{Perf}(\mathscr C)$ is $\mathscr I_{\mathscr C}$.
Since $[n]$ is an autoequivalence of $\mathsf{Perf}(\mathscr C)$, we obtain $\mathscr I_{\mathscr C}=\mathsf{Perf}(\mathscr C)$.
The Serre duality for $[n]$ is \eqref{eq:Nakayama-duality} restricted to perfect modules.
\end{proof}

We also recall the transfer from weak smooth to weak proper Calabi--Yau structures underlying \cite[Theorem 3.1]{Brav-Dyckerhoff2019}, which will be used in Theorem \ref{thm:koszul-CY-comparison}.

\begin{theorem}[{\cite[Theorem 3.1]{Brav-Dyckerhoff2019}}]\label{thm:smooth-proper-CY}
Let $\mathscr C$ be a smooth dg or strictly unital $A_\infty$ category.
A weak smooth $n$-Calabi--Yau structure on $\mathscr C$, in the sense of Definition \ref{def:weak-smooth-CY}, induces a weak proper $n$-Calabi--Yau structure on $\mathscr C^{\mathrm{prop}}$, and hence on every full subcategory of $\mathscr C^{\mathrm{prop}}$.	
\end{theorem}

For a Koszul-dualizing inclusion $j:\mathscr B\hookrightarrow\mathscr A$, condition (i) in Definition \ref{def:Koszul-duality} says that every $jb$ is right proper in $\mathscr A$.
A weak smooth $n$-Calabi--Yau structure on $\mathscr A$ in the sense of Definition \ref{def:weak-smooth-CY} therefore induces a weak proper $n$-Calabi--Yau structure on $\mathscr B$ by Theorem \ref{thm:smooth-proper-CY}.
Lemma \ref{lem:pCY-Gorenstein} then verifies the Gorenstein hypothesis in Proposition \ref{prop:restriction} and Theorem \ref{thm:main}.

\section{Koszul duality and reflexivity}\label{sec:koszul-reflexivity}

In this section, we prove that the functor $\Phi$ in \eqref{eq:compare-functor}, for a Koszul-dualizing inclusion $j:\mathscr B\hookrightarrow\mathscr A$, is quasi-fully faithful,
and furthermore that it is essentially surjective precisely when $\mathscr A$ is reflexive.
Under the Gorenstein condition \eqref{eq:Gorenstein-cond}, we also identify $\Phi$ with $\nu_{\mathscr B}\circ j^*$ and obtain the restriction and coinduction functors appearing in Theorem \ref{thm:main}.

\subsection{Koszul-dualizing inclusions}

Let $j:\mathscr B\hookrightarrow\mathscr A$ be a full inclusion.
The standard induction, restriction, and coinduction functors satisfy
\begin{equation}\label{eq:j-adjunction}
j_!\dashv j^*\dashv j_*;
\end{equation}
see \cite[Section 3.8]{Keller2006b}.
With the right-module convention of Section \ref{sec:basic}, these take the following explicit forms:
\begin{equation}\label{eq:j-module-formulas}
\begin{aligned}
(j^*M)(b)&=M(jb),\\
(j_!P)(a)&=P\otimes_{\mathscr B} \Gamma_j(a,-),\\
(j_*P)(a)&=\mathrm{Hom}_{\mathscr B}\bigl(j^*\mathbf Y_{\mathscr A}^{\mathrm r}(a),P\bigr),
\end{aligned}
\end{equation}
where $M\in\mathsf{Mod}(\mathscr A)$, $P\in\mathsf{Mod}(\mathscr B)$, and $\Gamma_j(a,b)=\mathrm{Hom}_{\mathscr A}(a,jb)$ as in the introduction.
Thus induction uses the left $\mathscr B$-module $\mathrm{Hom}_{\mathscr A}(a,j(-))$, whereas coinduction uses the right $\mathscr B$-module $\mathrm{Hom}_{\mathscr A}(j(-),a)$.
Fullness of $j$ implies that the unit $\mathrm{Id}\to j^*j_!$ and the counit $j^*j_*\to\mathrm{Id}$ are quasi-isomorphisms.
In particular, $j_!$ and $j_*$ are quasi-fully faithful.

We recall the notion of a Koszul-dualizing subcategory from \cite[Definition 2.12]{Ganatra-Gao-Venkatesh2022}, based on the hypotheses of \cite[Proposition 9.3]{Efimov2017}.

\begin{definition}\label{def:Koszul-duality}
Suppose that $\mathscr A$ is smooth and $\mathscr B$ is proper.
A full inclusion $j:\mathscr B\hookrightarrow\mathscr A$ is \textbf{Koszul-dualizing} if the following conditions hold:
\begin{enumerate}
    \item[(i)] Each $jb$ is right proper in $\mathscr A$, in the sense of Definition \ref{def:perf-prop}; equivalently, $\Gamma_j(a,b)\in\mathsf{Perf}(\Bbbk)$ for every $a\in\mathscr A$ and $b\in\mathscr B$.
    \item[(ii)] Restriction of left modules is quasi-fully faithful on perfect modules:
    \begin{equation}\label{eq:left-Koszul}
    (j^{\mathrm{op}})^*:\mathsf{Perf}(\mathscr A^{\mathrm{op}})\longrightarrow\mathsf{PsPerf}(\mathscr B^{\mathrm{op}}).
    \end{equation}
    \item[(iii)] For $P,Q\in\mathsf{PsPerf}(\mathscr B^{\mathrm{op}})$, evaluation of module morphisms induces a quasi-isomorphism
    \begin{equation}\label{eq:GGV-tensor}
    P^\vee\otimes_{\mathscr B} Q
    \longrightarrow\mathrm{Hom}_{\mathscr B^{\mathrm{op}}}(Q,P)^\vee.
    \end{equation}
\end{enumerate}
\end{definition}

Condition (i) ensures that the target in \eqref{eq:left-Koszul} is well defined.
For the same reason, the formula
\begin{equation}\label{eq:Phi-values}
\Phi(M)(b)\simeq\mathrm{Hom}_{\mathscr A}\bigl(M,\mathbf Y_{\mathscr A}^{\mathrm r}(jb)\bigr)^\vee,
\qquad M\in\mathsf{Perf}(\mathscr A),
\end{equation}
defines a pseudo-perfect right $\mathscr B$-module.
By the dualities \eqref{eq:perf-op-equivalence} and \eqref{eq:psperf-op-equivalence}, condition (ii) makes the functor $\Phi$ of \eqref{eq:compare-functor-Hom-form} quasi-fully faithful.

\begin{lemma}\label{lem:Koszul-generation}
Let $j:\mathscr B\hookrightarrow\mathscr A$ be Koszul-dualizing.
Induction and restriction give mutually inverse quasi-equivalences
\begin{equation}\label{eq:j_!-perfB-psperfA}
j_!:\mathsf{Perf}(\mathscr B)\xrightarrow{\;\sim\;}\mathsf{PsPerf}(\mathscr A),
\qquad
j^*:\mathsf{PsPerf}(\mathscr A)\xrightarrow{\;\sim\;}\mathsf{Perf}(\mathscr B).
\end{equation}
Consequently,
\begin{equation}\label{eq:Koszul-generation}
\mathsf{PsPerf}(\mathscr A)
=\mathrm{thick}_{\mathsf{Perf}(\mathscr A)}
\{\mathbf Y_{\mathscr A}^{\mathrm r}(jb):b\in\mathscr B\},
\end{equation}
and $\Phi$ restricts to a quasi-equivalence
\begin{equation}\label{eq:Phi-pseudo-perfect}
\Phi:\mathsf{PsPerf}(\mathscr A)\xrightarrow{\;\sim\;}\mathscr I_{\mathscr B},
\end{equation}
where $\mathscr I_{\mathscr B}$ is defined in \eqref{eq:Nakayama-image}.
\end{lemma}
\begin{proof}
By Definition \ref{def:Koszul-duality}(i), $\Gamma_j$ is pseudo-perfect, and fullness of $j$ makes $j_!$ quasi-fully faithful.
Conditions (ii)--(iii) give the remaining hypotheses of \cite[Proposition 9.3]{Efimov2017}, whose categorical proof applies over any field.
It follows that $j_!:\mathsf{Perf}(\mathscr B)\to\mathsf{PsPerf}(\mathscr A)$ is a quasi-equivalence, with quasi-inverse $j^*$ by \eqref{eq:j-adjunction}.
Since $j_!\mathbf Y_{\mathscr B}^{\mathrm r}(b)\simeq\mathbf Y_{\mathscr A}^{\mathrm r}(jb)$, the generation statement \eqref{eq:Koszul-generation} follows.

By fullness of $j$ and \eqref{eq:Phi-values},
\begin{equation}\label{eq:Phi-representables}
\Phi\bigl(\mathbf Y_{\mathscr A}^{\mathrm r}(jb)\bigr)
\simeq\mathbf Y_{\mathscr B}^{\mathrm l}(b)^\vee
\simeq\nu_{\mathscr B}\mathbf Y_{\mathscr B}^{\mathrm r}(b).
\end{equation}
Since $\Phi$ is quasi-fully faithful, \eqref{eq:Koszul-generation} identifies its essential image on $\mathsf{PsPerf}(\mathscr A)$ with the thick subcategory generated by these modules, which is $\mathscr I_{\mathscr B}$ as defined in \eqref{eq:Nakayama-image}.
\end{proof}

\subsection{Relation to the Nakayama functor}

The Gorenstein condition allows us to recover restriction from $\Phi$ by the inverse Nakayama functor.
This identifies the additional functors in Theorem \ref{thm:main} without changing the essential-surjectivity question.

\begin{proposition}\label{prop:restriction}
Let $j:\mathscr B\hookrightarrow\mathscr A$ be Koszul-dualizing, and suppose that $\mathscr B$ satisfies \eqref{eq:Gorenstein-cond}.
There is a natural quasi-isomorphism
\begin{equation}\label{eq:Phi-Nakayama-restriction}
\Phi\simeq\nu_{\mathscr B}\circ j^*
\qquad\text{on }\mathsf{Perf}(\mathscr A).
\end{equation}
In particular, restriction gives a quasi-fully faithful functor
\[
j^*:\mathsf{Perf}(\mathscr A)\longrightarrow\mathsf{PsPerf}(\mathscr B),
\]
whose restriction to $\mathsf{PsPerf}(\mathscr A)$ is the quasi-inverse of $j_!$ in \eqref{eq:j_!-perfB-psperfA}.
The functor $\Phi$ is essentially surjective if and only if coinduction restricts to a quasi-equivalence
\begin{equation}\label{eq:j_*PsPerfB-in-PerfA}
j_*:\mathsf{PsPerf}(\mathscr B)\xrightarrow{\;\sim\;}\mathsf{Perf}(\mathscr A).
\end{equation}
When these conditions hold, $j^*$ and $j_*$ are mutually inverse on these subcategories.
\end{proposition}
\begin{proof}
For $M\in\mathsf{Perf}(\mathscr A)$ and $b\in\mathscr B$, quasi-full faithfulness of $\Phi$, \eqref{eq:Phi-representables}, and the adjunction $\nu_{\mathscr B}\dashv\nu_{\mathscr B}^{\mathsf R}$ give
\begin{align*}
(j^*M)(b)
&\simeq\mathrm{Hom}_{\mathscr A}\bigl(\mathbf Y_{\mathscr A}^{\mathrm r}(jb),M\bigr)\\
&\simeq\mathrm{Hom}_{\mathscr B}\bigl(\nu_{\mathscr B}\mathbf Y_{\mathscr B}^{\mathrm r}(b),\Phi(M)\bigr)\\
&\simeq\mathrm{Hom}_{\mathscr B}\bigl(\mathbf Y_{\mathscr B}^{\mathrm r}(b),\nu_{\mathscr B}^{\mathsf R}\Phi(M)\bigr)\\
&\simeq\bigl(\nu_{\mathscr B}^{\mathsf R}\Phi(M)\bigr)(b).
\end{align*}
These identifications are natural in $M$ and $b$, so $j^*\simeq\nu_{\mathscr B}^{\mathsf R}\circ\Phi$.
By Lemma \ref{lem:Nakayama-Serre}, the Gorenstein condition makes $\nu_{\mathscr B}$ an autoequivalence of $\mathsf{Mod}(\mathscr B)$ with inverse $\nu_{\mathscr B}^{\mathsf R}$, and both functors preserve $\mathsf{PsPerf}(\mathscr B)$.
This proves \eqref{eq:Phi-Nakayama-restriction} and quasi-full faithfulness of $j^*$ on $\mathsf{Perf}(\mathscr A)$.
Its assertion on pseudo-perfect modules is \eqref{eq:j_!-perfB-psperfA}.

Since $\nu_{\mathscr B}$ restricts to an autoequivalence of $\mathsf{PsPerf}(\mathscr B)$, \eqref{eq:Phi-Nakayama-restriction} shows that $\Phi$ is essentially surjective precisely when $j^*$ is.
For $M\in\mathsf{Perf}(\mathscr A)$, the unit $M\to j_*j^*M$ is a quasi-isomorphism: at $a\in\mathscr A$, this follows from
\[
(j_*j^*M)(a)
\simeq\mathrm{Hom}_{\mathscr B}\bigl(j^*\mathbf Y_{\mathscr A}^{\mathrm r}(a),j^*M\bigr)
\simeq\mathrm{Hom}_{\mathscr A}\bigl(\mathbf Y_{\mathscr A}^{\mathrm r}(a),M\bigr)
\simeq M(a).
\]
If $j^*$ is essentially surjective, every $P\in\mathsf{PsPerf}(\mathscr B)$ has the form $j^*M$ up to quasi-isomorphism, and hence $j_*P\simeq M$ is perfect.
The adjunction $j^*\dashv j_*$ therefore restricts to the quasi-equivalences asserted in \eqref{eq:j_*PsPerfB-in-PerfA}.
Conversely, if \eqref{eq:j_*PsPerfB-in-PerfA} holds, then $j^*j_*P\simeq P$ for every pseudo-perfect $P$ by fullness of $j$, proving essential surjectivity of $j^*$ and of $\Phi$.
\end{proof}

\subsection{Reflexivity}

To characterize essential surjectivity of $\Phi$, we relate it to the notion of reflexivity in the sense of \cite[Definition 2.1.1 and Remark 2.1.5(2)]{Booth-Goodbody-Opper2025}.
For a small dg or strictly unital $A_\infty$ category $\mathscr C$, Lemma \ref{lem:module-finiteness} ensures that the evaluation functor
\begin{equation}\label{eq:evaluation}
\begin{aligned}
\operatorname{ev}_{\mathscr C}:\mathsf{Perf}(\mathscr C)^{\mathrm{op}}
&\longrightarrow\mathsf{PsPerf}\bigl(\mathsf{PsPerf}(\mathscr C)^{\mathrm{op}}\bigr),\\
M&\longmapsto\bigl[N\longmapsto\mathrm{Hom}_{\mathscr C}(M,N)\bigr]
\end{aligned}
\end{equation}
is well defined. 

\begin{definition}\label{def:reflexivity}
The category $\mathscr C$ is \textbf{reflexive} if \eqref{eq:evaluation} is a quasi-equivalence, and \textbf{semireflexive} if it is quasi-fully faithful.
\end{definition}

Reflexivity is Morita invariant \cite[Proposition 3.8]{Goodbody2026} and invariant under passage to the opposite category \cite[Remark 2.1.5(1)]{Booth-Goodbody-Opper2025}.
For $A_\infty$ categories, we use the dg models fixed in Section \ref{sec:basic}; the induced equivalences of module categories identify their evaluation functors.

\begin{proposition}\label{prop:reflexivity}
Let $j:\mathscr B\hookrightarrow\mathscr A$ be Koszul-dualizing.
Then $\mathscr A$ is semireflexive, and the functor $\Phi$ in \eqref{eq:compare-functor} is essentially surjective if and only if $\mathscr A$ is reflexive.
\end{proposition}
\begin{proof}
By Lemma \ref{lem:Koszul-generation}, the functor
\[
i:\mathscr B\longrightarrow\mathsf{PsPerf}(\mathscr A),
\qquad b\longmapsto\mathbf Y_{\mathscr A}^{\mathrm r}(jb),
\]
is a Morita equivalence.
Restriction along $i^{\mathrm{op}}$ therefore induces a quasi-equivalence
\begin{equation}\label{eq:reflexivity-Morita}
(i^{\mathrm{op}})^*:
\mathsf{PsPerf}\bigl(\mathsf{PsPerf}(\mathscr A)^{\mathrm{op}}\bigr)
\xrightarrow{\;\sim\;}\mathsf{PsPerf}(\mathscr B^{\mathrm{op}});
\end{equation}
see \cite[Proposition 3.5]{Goodbody2026}.
For $M\in\mathsf{Perf}(\mathscr A)$, the composite $(i^{\mathrm{op}})^*\operatorname{ev}_{\mathscr A}$ takes the value
\[
\bigl((i^{\mathrm{op}})^*\operatorname{ev}_{\mathscr A}(M)\bigr)(b)
=\mathrm{Hom}_{\mathscr A}\bigl(M,\mathbf Y_{\mathscr A}^{\mathrm r}(jb)\bigr).
\]
It is thus naturally identified with
\begin{equation}\label{eq:mixed-Hom-functor}
\mathsf{Perf}(\mathscr A)^{\mathrm{op}}
\xrightarrow[\sim]{\mathrm{Hom}_{\mathscr A}(-,\mathscr A)}
\mathsf{Perf}(\mathscr A^{\mathrm{op}})
\xrightarrow{(j^{\mathrm{op}})^*}
\mathsf{PsPerf}(\mathscr B^{\mathrm{op}}).
\end{equation}
The first arrow is the duality \eqref{eq:perf-op-equivalence}, and the second is quasi-fully faithful by Definition \ref{def:Koszul-duality}(ii).
Together with \eqref{eq:reflexivity-Morita}, this proves that $\operatorname{ev}_{\mathscr A}$ is quasi-fully faithful.
Moreover, $\operatorname{ev}_{\mathscr A}$ is essentially surjective precisely when \eqref{eq:mixed-Hom-functor} is.
Taking the opposite of \eqref{eq:mixed-Hom-functor} and applying the linear duality \eqref{eq:psperf-op-equivalence} gives $\Phi$, as in \eqref{eq:compare-functor-Hom-form}.
These dualities preserve essential surjectivity, proving the assertion.
\end{proof}

\section{Boundary categories and Koszul duality}\label{sec:infinity-singularity}

We describe the two quotient categories in Theorem \ref{thm:main}, then use the results of Section \ref{sec:koszul-reflexivity} to construct the comparison functor and characterize when it is a quasi-equivalence.
We also prove Theorem \ref{thm:koszul-CY-comparison} by transporting the categorical residue through this equivalence.
The final subsection recalls a trace on the singularity category obtained from Hochschild localization.

\subsection{Categorical formal punctured neighborhood of infinity}\label{subsec:formal-punctured-neighborhood}

The $\Bbbk$-linear Calkin dg category is the Drinfeld quotient
\[
q:\mathsf{Mod}(\Bbbk)\longrightarrow\mathsf{Calk}(\Bbbk)
:=\mathsf{Mod}(\Bbbk)/\mathsf{Perf}(\Bbbk);
\]
see \cite[Proposition 1.1]{Efimov2017} or \cite[Section 2.2]{Ganatra-Gao-Venkatesh2022}.
For a small dg or strictly unital $A_\infty$ category $\mathscr A$, composition with $q$ gives the Calkin-valued Yoneda functor
\begin{equation}\label{eq:reduced-Yoneda}
\begin{aligned}
\overline{\mathbf Y}:\mathscr A
&\longrightarrow\mathsf{Fun}(\mathscr A^{\mathrm{op}},\mathsf{Calk}(\Bbbk)),\\
a&\longmapsto q\circ\mathbf Y_{\mathscr A}^{\mathrm r}(a).
\end{aligned}
\end{equation}
Here $\mathsf{Fun}$ denotes the dg category of strictly unital $A_\infty$ functors and $A_\infty$ natural transformations, which models the homotopy internal Hom; see \cite[Section 1 and Proposition 3.6]{Efimov2017} and \cite[Sections 2.1 and 2.3]{Ganatra-Gao-Venkatesh2022}.

\begin{definition}\label{def:cat-form-punc-nbhd}
The \textbf{categorical formal punctured neighborhood of infinity} of $\mathscr A$ is the full essential image
\[
\widehat{\mathscr A}_{\infty}
:=\operatorname{ess\,im}(\overline{\mathbf Y})
\subset\mathsf{Fun}(\mathscr A^{\mathrm{op}},\mathsf{Calk}(\Bbbk)).
\]
We use the full essential image before pretriangulated or idempotent completion, as in \cite[Definition 3.3(3)]{Efimov2017}.
\end{definition}

Since $q$ annihilates precisely the perfect $\Bbbk$-complexes, an object $a\in\mathscr A$ satisfies
\[
\overline{\mathbf Y}(a)\simeq0
\quad\Longleftrightarrow\quad
a\in\mathscr A^{\mathrm{prop}},
\]
where right-proper objects are defined in Definition \ref{def:perf-prop}; see also \cite[Lemma 2.8]{Ganatra-Gao-Venkatesh2022}.
When $\mathscr A$ is smooth, the right-proper objects of $\mathsf{Perf}(\mathscr A)$ are exactly $\mathsf{PsPerf}(\mathscr A)$ by Lemma \ref{lem:module-finiteness}; compare \cite[Proposition 4.1]{Efimov2017}.
Under the Koszul-duality hypotheses, the functor induced on the corresponding quotient is a quasi-equivalence.

\begin{proposition}[{\cite[Proposition 9.3]{Efimov2017},\cite[Proposition 2.13]{Ganatra-Gao-Venkatesh2022}}] \label{prop:boundary-quotient}
Let $\mathscr A$ be smooth and let $j:\mathscr B\hookrightarrow\mathscr A$ be Koszul-dualizing in the sense of Definition \ref{def:Koszul-duality}.
The Calkin-valued Yoneda functor induces quasi-equivalences
\begin{equation}\label{eq:boundary-quotient}
\mathscr A/\mathscr B\xrightarrow{\;\sim\;}\widehat{\mathscr A}_{\infty}
\end{equation}
and
\begin{equation}\label{eq:perf-boundary-quotient}
\mathsf{Perf}(\mathscr A)/\mathsf{PsPerf}(\mathscr A)
\xrightarrow{\;\sim\;}\widehat{\mathsf{Perf}(\mathscr A)}_{\infty}.
\end{equation}
The induced functor $\widehat{\mathscr A}_{\infty}\to\widehat{\mathsf{Perf}(\mathscr A)}_{\infty}$ is a Morita equivalence.
\end{proposition}

\subsection{Singularity category}\label{subsec:singularity}

Let $\mathscr B$ be a proper dg or strictly unital $A_\infty$ category, so that $\mathsf{Perf}(\mathscr B)\subset\mathsf{PsPerf}(\mathscr B)$ by Lemma \ref{lem:module-finiteness}.
We use the dg version of the singularity category recalled in \cite[Section 2.1]{Efimov2017}.
In that reference, local properness is the properness condition of Definition \ref{def:smooth-proper}.

\begin{definition}\label{def:dg-singularity}
The \textbf{dg singularity category} of $\mathscr B$ is the Drinfeld quotient
\begin{equation}\label{eq:dg-singularity}
\mathscr D_{\mathsf{sg}}(\mathscr B)
:=\mathsf{PsPerf}(\mathscr B)/\mathsf{Perf}(\mathscr B).
\end{equation}
\end{definition}

Its homotopy category is the Verdier quotient
\[
H^0(\mathscr D_{\mathsf{sg}}(\mathscr B))
\simeq H^0(\mathsf{PsPerf}(\mathscr B))/H^0(\mathsf{Perf}(\mathscr B))
\]
by \cite[Theorem 3.4]{Drinfeld2004}.
Since $\mathsf{Perf}(\mathscr B)$ is thick in $\mathsf{PsPerf}(\mathscr B)$, the quotient vanishes precisely when every pseudo-perfect $\mathscr B$-module is perfect.
In particular, it vanishes when $\mathscr B$ is smooth, by Lemma \ref{lem:module-finiteness}.

For a Koszul-dualizing inclusion $j:\mathscr B\hookrightarrow\mathscr A$, the comparison functor $\Phi$ sends $\mathsf{PsPerf}(\mathscr A)$ to $\mathscr I_{\mathscr B}$ by \eqref{eq:Phi-pseudo-perfect}.
The following consequence of linear duality identifies the quotient by this subcategory.

\begin{lemma}\label{lem:Dsg-op-op}
Linear duality induces a quasi-equivalence
\begin{equation}\label{eq:singularity-duality}
\mathsf{PsPerf}(\mathscr B)/\mathscr I_{\mathscr B}
\xrightarrow{\;\sim\;}
\mathscr D_{\mathsf{sg}}(\mathscr B^{\mathrm{op}})^{\mathrm{op}},
\qquad [N]\longmapsto[N^\vee],
\end{equation}
where $\mathscr I_{\mathscr B}$ is defined in \eqref{eq:Nakayama-image}.
\end{lemma}
\begin{proof}
The duality \eqref{eq:psperf-op-equivalence} identifies $\mathsf{PsPerf}(\mathscr B)$ with $\mathsf{PsPerf}(\mathscr B^{\mathrm{op}})^{\mathrm{op}}$.
It sends each generator $\mathbf Y_{\mathscr B}^{\mathrm l}(b)^\vee$ of $\mathscr I_{\mathscr B}$ to the representable left module $\mathbf Y_{\mathscr B}^{\mathrm l}(b)$, since $\mathscr B$ is proper.
It therefore identifies $\mathscr I_{\mathscr B}$ with $\mathsf{Perf}(\mathscr B^{\mathrm{op}})^{\mathrm{op}}$ and induces \eqref{eq:singularity-duality} on dg quotients.
Pointwise double duality gives the inverse.
\end{proof}

If $\mathscr B$ satisfies the Gorenstein condition \eqref{eq:Gorenstein-cond}, then $\mathscr I_{\mathscr B}=\mathsf{Perf}(\mathscr B)$, so Lemma \ref{lem:Dsg-op-op} gives a linear-duality quasi-equivalence $\mathscr D_{\mathsf{sg}}(\mathscr B^{\mathrm{op}})^{\mathrm{op}}\simeq\mathscr D_{\mathsf{sg}}(\mathscr B)$.
For a weak proper $n$-Calabi--Yau structure, Lemma \ref{lem:pCY-Gorenstein} further gives $\nu_{\mathscr B}^{-1}\simeq[-n]$.
Composing linear duality with this inverse Nakayama functor gives the normalization
\begin{equation}\label{eq:CY-singularity-duality}
(-)^\vee[-n]:
\mathscr D_{\mathsf{sg}}(\mathscr B^{\mathrm{op}})^{\mathrm{op}}
\xrightarrow{\;\sim\;}\mathscr D_{\mathsf{sg}}(\mathscr B).
\end{equation}

\subsection{Comparison under Koszul duality}\label{subsec:koszul-comparison}

We now apply $\Phi$ to the quotient presentations \eqref{eq:perf-boundary-quotient} and \eqref{eq:singularity-duality} to prove Theorem \ref{thm:main}.

\begin{proof}[Proof of Theorem \ref{thm:main}]
By Definition \ref{def:Koszul-duality}(ii) and the factorization \eqref{eq:compare-functor-Hom-form}, the functor $\Phi$ is quasi-fully faithful.
Lemma \ref{lem:Koszul-generation} identifies its restriction to pseudo-perfect modules with a quasi-equivalence
\[
\Phi:\mathsf{PsPerf}(\mathscr A)\xrightarrow{\;\sim\;}\mathscr I_{\mathscr B}.
\]
Thus, after identifying $\mathsf{Perf}(\mathscr A)$ with the full essential image of $\Phi$, we have full pretriangulated subcategories
\[
\mathscr I_{\mathscr B}
\subset\operatorname{ess\,im}(\Phi)
\subset\mathsf{PsPerf}(\mathscr B).
\]
The essential image is thick on homotopy categories: exactness and full faithfulness give closure under shifts and cones, while idempotents lift to $H^0(\mathsf{Perf}(\mathscr A))$ and split there.
Lemma \ref{lem:quotient-density} therefore gives a quasi-fully faithful functor on the corresponding dg quotients.
Using Proposition \ref{prop:boundary-quotient} and Lemma \ref{lem:Dsg-op-op}, we obtain
\begin{equation}\label{eq:comparison-quotients}
\begin{aligned}
\widehat{\mathsf{Perf}(\mathscr A)}_{\infty}
&\simeq\mathsf{Perf}(\mathscr A)/\mathsf{PsPerf}(\mathscr A)\\
&\xrightarrow{\;\overline\Phi\;}
\mathsf{PsPerf}(\mathscr B)/\mathscr I_{\mathscr B}
\xrightarrow[\sim]{\;(-)^\vee\;}
\mathscr D_{\mathsf{sg}}(\mathscr B^{\mathrm{op}})^{\mathrm{op}},
\end{aligned}
\end{equation}
where $\overline\Phi$ is induced by $\Phi$.
The composite is the functor \eqref{eq:main-comparison}.

The same application of Lemma \ref{lem:quotient-density} shows that $\overline\Phi$ is essentially surjective if and only if $\Phi$ is, proving \textup{(i)} $\Longleftrightarrow$ \textup{(ii)}.
Proposition \ref{prop:reflexivity} identifies essential surjectivity of $\Phi$ with reflexivity of $\mathscr A$, giving \textup{(ii)} $\Longleftrightarrow$ \textup{(iii)}.

Suppose now that $\mathscr B$ satisfies the Gorenstein condition \eqref{eq:Gorenstein-cond}.
Then $\mathscr I_{\mathscr B}=\mathsf{Perf}(\mathscr B)$, and Lemma \ref{lem:Nakayama-Serre} makes $\nu_{\mathscr B}$ an autoequivalence of $\mathsf{Mod}(\mathscr B)$ preserving both $\mathsf{Perf}(\mathscr B)$ and $\mathsf{PsPerf}(\mathscr B)$.
Proposition \ref{prop:restriction} gives
\[
\Phi\simeq\nu_{\mathscr B}\circ j^*,
\qquad
j^*\simeq\nu_{\mathscr B}^{-1}\circ\Phi
\qquad\text{on }\mathsf{Perf}(\mathscr A).
\]
Hence $j^*$ is quasi-fully faithful, and \eqref{eq:j_!-perfB-psperfA} identifies its restriction to $\mathsf{PsPerf}(\mathscr A)$ with a quasi-equivalence onto $\mathsf{Perf}(\mathscr B)$.
On the dg quotients, the functor induced by $j^*$ is therefore $\overline\Phi$ followed by the autoequivalence induced by $\nu_{\mathscr B}^{-1}$.
This proves the assertions about \eqref{eq:ordinary-comparison}.
Finally, Proposition \ref{prop:restriction} identifies condition \textup{(ii)} with the assertion that $j_*$ restricts to the quasi-equivalence \eqref{eq:j_*PsPerfB-in-PerfA}, with quasi-inverse $j^*$.
\end{proof}

We next use the restriction functor of Theorem \ref{thm:main} to transport the categorical residue of \cite[Theorem 1.7]{Gao2025} to the singularity category.

\begin{theorem}\label{thm:koszul-CY-comparison}
Suppose that $\mathscr A$ carries a weak smooth $n$-Calabi--Yau structure and is reflexive, and let $j:\mathscr B\hookrightarrow\mathscr A$ be a Koszul-dualizing full inclusion.
Restriction induces a quasi-equivalence
\begin{equation}\label{eq:CY-boundary-comparison}
\overline{j^*}:\widehat{\mathsf{Perf}(\mathscr A)}_{\infty}
\xrightarrow{\;\sim\;}\mathscr D_{\mathsf{sg}}(\mathscr B).
\end{equation}
The categorical residue transports through this quasi-equivalence to a weak generalized proper $(n-1)$-Calabi--Yau structure on $\mathscr D_{\mathsf{sg}}(\mathscr B)$, in the sense of Remark \ref{rem:CY-finiteness}.
\end{theorem}
\begin{proof}
Condition (i) in Definition \ref{def:Koszul-duality} places $j\mathscr B$ in $\mathscr A^{\mathrm{prop}}$.
Theorem \ref{thm:smooth-proper-CY} therefore gives a weak proper $n$-Calabi--Yau structure on $\mathscr B$.
By Lemma \ref{lem:pCY-Gorenstein}, the category $\mathscr B$ satisfies the Gorenstein condition and $\nu_{\mathscr B}\simeq[n]$.
Since $\mathscr A$ is reflexive, Theorem \ref{thm:main} gives the quasi-equivalence \eqref{eq:CY-boundary-comparison}.
Moreover, $\Phi\simeq[n]\circ j^*$ by \eqref{eq:Phi-Nakayama-restriction}, so composing \eqref{eq:main-comparison} with the duality \eqref{eq:CY-singularity-duality} gives this same restriction functor: double duality recovers $\Phi$, and the shift $[-n]$ cancels its Nakayama shift.

The weak smooth structure extends to $\mathsf{Perf}(\mathscr A)$ by Morita invariance.
By \cite[Theorem 1.7 and equation (3.67)]{Gao2025}, it determines a residue trace
\begin{equation}\label{eq:boundary-residue}
\operatorname{res}:C_*(\widehat{\mathsf{Perf}(\mathscr A)}_{\infty})
\longrightarrow\Bbbk[1-n]
\end{equation}
whose pairings on morphism cohomology are nondegenerate.
The quasi-isomorphism $C_*(\overline{j^*})$ identifies homotopy classes of traces by precomposition.
There is therefore a unique homotopy class of traces on $\mathscr D_{\mathsf{sg}}(\mathscr B)$ whose pullback is the class of \eqref{eq:boundary-residue}.
Its pairings are nondegenerate because $\overline{j^*}$ identifies morphism cohomology and composition.
This gives the asserted weak generalized proper Calabi--Yau structure.
\end{proof}

For a proper scheme $X$ over a perfect field $\Bbbk$, let $\mathsf{Coh}(X)$ denote the dg category of complexes with bounded coherent cohomology and let $\omega_X$ be its dualizing complex.
Efimov's geometric comparison is
\begin{equation}\label{eq:scheme-boundary-comparison}
\widehat{\mathsf{Coh}(X)}_{\infty}
\simeq
\frac{\mathsf{Coh}(X)}{\omega_X\otimes_{\mathcal O_X} \mathsf{Perf}(X)}
\simeq\mathscr D_{\mathsf{sg}}(X)^{\mathrm{op}};
\end{equation}
see \cite[Theorem 9.2 and Proposition B.1]{Efimov2017} and \cite[Section 4]{Gao2025}.
The denominator is the full essential image of tensoring perfect complexes with $\omega_X$, and the second quasi-equivalence is induced by Grothendieck--Verdier duality.
For Gorenstein $X$, the denominator is $\mathsf{Perf}(X)$; if moreover $X$ has dimension $n$ and $\omega_X\simeq\mathcal O_X[n]$, then \cite[Theorem 1.6]{Gao2025} gives the corresponding weak generalized proper $(n-1)$-Calabi--Yau structure on $\mathscr D_{\mathsf{sg}}(X)$.
In characteristic zero, the chosen trivialization supplies the weak smooth structure on $\mathsf{Coh}(X)$ by \cite[Section 5.2 and Proposition 5.12]{Brav-Dyckerhoff2019}.

\subsection{Traces on the singularity category}\label{subsec:singularity-traces}

For a proper category $\mathscr B$ with a weak proper $n$-Calabi--Yau structure, Hochschild localization gives a trace directly on $\mathscr D_{\mathsf{sg}}(\mathscr B)$.
We recall this construction and the finiteness condition under which its pairings are perfect.
By Morita invariance, regard the trace \eqref{eq:pCY-trace} as a trace on $\mathsf{Perf}(\mathscr B)$.
Localization for Hochschild chains gives an exact triangle
\begin{equation}\label{eq:singularity-Hochschild-triangle}
\begin{aligned}
C_*(\mathsf{Perf}(\mathscr B))&\longrightarrow C_*(\mathsf{PsPerf}(\mathscr B))
\longrightarrow C_*(\mathscr D_{\mathsf{sg}}(\mathscr B))\\
&\xrightarrow{\;\partial\;}C_*(\mathsf{Perf}(\mathscr B))[1];
\end{aligned}
\end{equation}
see \cite[Theorem 5.2(a)]{Keller2006b}.
Thus the induced trace is the composite
\begin{equation}\label{eq:singularity-trace}
\operatorname{tr}_{\mathsf{sg}}
:=\operatorname{tr}[1]\circ\partial:
C_*(\mathscr D_{\mathsf{sg}}(\mathscr B))\longrightarrow\Bbbk[1-n].
\end{equation}
The following nondegeneracy statement is a specialization of \cite[Theorem 4.3.1(b)]{Hanihara-Liu2026}.

\begin{theorem}\label{thm:singularity-CY}
Let $\mathscr B$ carry a weak proper $n$-Calabi--Yau structure, and suppose that
\begin{equation}\label{eq:PsPerf-local-finiteness}
\dim_{\Bbbk}H^i\mathrm{Hom}_{\mathscr B}(M,N)<\infty
\qquad
\text{for all }M,N\in\mathsf{PsPerf}(\mathscr B)\text{ and }i\in\mathbb Z.
\end{equation}
Then the trace \eqref{eq:singularity-trace} induces a bimodule quasi-isomorphism
\[
\mathscr D_{\mathsf{sg}}(\mathscr B)_\Delta
\xrightarrow{\;\sim\;}\mathscr D_{\mathsf{sg}}(\mathscr B)^\vee[1-n].
\]
It therefore defines a weak locally proper $(n-1)$-Calabi--Yau structure in the sense of Remark \ref{rem:CY-finiteness}.
\end{theorem}
\begin{proof}
Put $\mathscr C=\mathsf{PsPerf}(\mathscr B)$ and $\mathscr P=\mathsf{Perf}(\mathscr B)$.
Condition \eqref{eq:PsPerf-local-finiteness} makes the canonical map $\mathscr C_\Delta\to\mathscr C_\Delta^{\vee\vee}$ a quasi-isomorphism.
Moreover, the restricted Yoneda functor
\[
\mathscr C\longrightarrow\mathsf{Mod}(\mathscr P),
\qquad M\longmapsto\bigl[P\longmapsto\mathrm{Hom}_{\mathscr B}(P,M)\bigr],
\]
is quasi-fully faithful, since it is the restriction of the quasi-equivalence of dg module categories induced by $\mathscr B\to\mathsf{Perf}(\mathscr B)$.
These verify the two hypotheses of \cite[Theorem 4.3.1(b)]{Hanihara-Liu2026} for $\mathscr P\subset\mathscr C$.
Apply that theorem with the shifted diagonal $\mathscr C_\Delta[n]$, corresponding to the shift autoequivalence $[-n]$ in its notation.
The weak proper Calabi--Yau structure gives $\mathscr P_\Delta[n]\simeq\mathscr P^\vee$, and the theorem identifies its image under the connecting morphism with a quasi-isomorphism
\[
\mathscr D_{\mathsf{sg}}(\mathscr B)_\Delta[n]
\xrightarrow{\;\sim\;}\mathscr D_{\mathsf{sg}}(\mathscr B)^\vee[1].
\]
Shifting by $[-n]$ gives the assertion.
\end{proof}

The finiteness hypothesis \eqref{eq:PsPerf-local-finiteness} cannot be omitted from this general statement.

\begin{example}\label{ex:singularity-CY-finiteness}
Let $\mathscr B$ have objects $b_i$, $i\geqslant1$, with $\mathrm{Hom}_{\mathscr B}(b_i,b_j)=\Bbbk\delta_{ij}$ concentrated in degree zero.
It is proper, and the trace
\[
C_*(\mathscr B)\simeq\bigoplus_{i\geqslant1}\Bbbk\longrightarrow\Bbbk,
\qquad (\lambda_i)_i\longmapsto\sum_i\lambda_i,
\]
defines a weak proper $0$-Calabi--Yau structure.
A pseudo-perfect module is a sequence of perfect $\Bbbk$-complexes; it is perfect over $\mathscr B$ precisely when only finitely many components are nonzero up to quasi-isomorphism.
Thus $P=(\Bbbk,\Bbbk,\ldots)$ is pseudo-perfect, and \eqref{eq:quotient-Hom} gives
\[
\mathrm{Hom}_{\mathscr D_{\mathsf{sg}}(\mathscr B)}(P,P)
\simeq
\operatorname{cone}\left(\bigoplus_{i\geqslant1}\Bbbk
\hookrightarrow\prod_{i\geqslant1}\Bbbk\right).
\]
Its cohomology is concentrated in degree zero, where it is the nonzero vector space $\prod_i\Bbbk/\bigoplus_i\Bbbk$.
A weak generalized proper $(-1)$-Calabi--Yau structure would give a nondegenerate pairing between $H^{0}\mathrm{Hom}_{\mathscr D_{\mathsf{sg}}(\mathscr B)}(P,P)\neq 0$ and $H^{-1}\mathrm{Hom}_{\mathscr D_{\mathsf{sg}}(\mathscr B)}(P,P)=0$.
Hence no such structure exists.
\end{example}

\section{Algebraic criteria}\label{sec:algebraic-criteria}

We express the hypotheses of Theorem \ref{thm:main} in terms of endomorphism algebras, prove Corollary \ref{cor:main}, and apply the comparison to Koszul duals of proper algebras.
We use the dg module models fixed in Section \ref{sec:basic}, and write $\mathrm{Hom}$ for their morphism complexes.
For an algebra $B$, the category $\mathsf{PsPerf}(B)$ consists of modules with finite-dimensional total cohomology.

\subsection{Split generators and endomorphism algebras}

Let $j:\mathscr{B}\hookrightarrow\mathscr{A}$ be a full inclusion, with $\mathscr{A}$ smooth and $\mathscr{B}$ proper. Choose split-generators $G\in\mathsf{Perf}(\mathscr{A})$ and $K\in\mathsf{Perf}(\mathscr{B})$, and put
\begin{equation}\label{eq:algebraic-models}
A=\mathrm{Hom}_{\mathscr{A}}(G,G), \qquad
B=\mathrm{Hom}_{\mathscr{B}}(K,K),\qquad
M=\mathrm{Hom}_{\mathscr{B}}(K,\Phi(G)).
\end{equation}
Smoothness supplies $G$ by \cite[Lemma 2.5]{Efimov2010}; the existence of $K$ is an additional assumption.
The Morita equivalences induced by $G$ and $K$ identify $\Phi$ with
\[
-\otimes_A M:\mathsf{Perf}(A)\longrightarrow\mathsf{PsPerf}(B),
\]
where $M$ is an $A$--$B$-bimodule with action $\alpha m\beta=\Phi(\alpha)\circ m\circ\beta$; see \cite[Section 3]{Keller2006b}.

\begin{proposition}\label{prop:algebraic-comparison}
Assume that each $jb$ is right proper in $\mathscr A$, as in Definition \ref{def:Koszul-duality}(i).
Then $\Phi$ takes values in $\mathsf{PsPerf}(\mathscr B)$, and it is quasi-fully faithful if and only if its action morphism
\begin{equation}\label{eq:algebraic-action}
A\longrightarrow\mathrm{Hom}_B(M,M),\qquad
\alpha\longmapsto\bigl[m\longmapsto\Phi(\alpha)\circ m\bigr],
\end{equation}
is a quasi-isomorphism.
In this case its essential image, under the preceding Morita equivalence, is $\mathrm{thick}_B(M)$.

Suppose that \eqref{eq:algebraic-action} is a quasi-isomorphism.
If
\begin{equation}\label{eq:algebraic-density-test}
\mathsf{PsPerf}(B)=\mathrm{thick}_B(M),
\end{equation}
then $\Phi$ is essentially surjective and hence a quasi-equivalence.
If, moreover,
\begin{equation}\label{eq:algebraic-degreewise-finiteness}
\dim_{\Bbbk}H^d(A)<\infty
\qquad\text{for every }d\in\mathbb Z,
\end{equation}
then $j$ is Koszul-dualizing and $\mathscr A$ is reflexive.
\end{proposition}
\begin{proof}
The right-properness assumption and \eqref{eq:Phi-values} imply that $\Phi$ takes values in $\mathsf{PsPerf}(\mathscr B)$.
The map \eqref{eq:algebraic-action} is the map on endomorphism complexes of the generator $A$.
Closure under shifts, finite cones, and direct summands in both variables shows that $\Phi$ is quasi-fully faithful precisely when this map is a quasi-isomorphism.
In that case its essential image is $\mathrm{thick}_B(M)$.
By \eqref{eq:compare-functor-Hom-form}, quasi-full faithfulness of $\Phi$ is equivalent to condition (ii) of Definition \ref{def:Koszul-duality}.

Assume that \eqref{eq:algebraic-action} is a quasi-isomorphism and that \eqref{eq:algebraic-density-test} holds.
The description of the essential image then shows that $\Phi$ is essentially surjective.
Thus it induces a quasi-equivalence
\[
\mathsf{Perf}(A)\xrightarrow{\;\sim\;}\mathsf{PsPerf}(B).
\]

Assume in addition \eqref{eq:algebraic-degreewise-finiteness}.
Closure under shifts, finite cones, and direct summands extends degreewise finiteness from $A$ to all morphism complexes in $\mathsf{Perf}(A)$.
The preceding quasi-equivalence and linear duality therefore imply this finiteness for all morphism complexes in $\mathsf{PsPerf}(B^{\mathrm{op}})$.
For $P,Q\in\mathsf{PsPerf}(B^{\mathrm{op}})$, tensor--Hom adjunction gives
\begin{equation}\label{eq:algebraic-biduality}
\bigl(P^\vee\otimes_B Q\bigr)^\vee
\simeq\mathrm{Hom}_{B^{\mathrm{op}}}(Q,P).
\end{equation}
Under this identification, the map \eqref{eq:GGV-tensor} is the canonical biduality map for $P^\vee\otimes_B Q$.
Degreewise finiteness makes this map a quasi-isomorphism, proving condition (iii).
Condition (i) holds by assumption, and condition (ii) follows from quasi-full faithfulness of $\Phi$.
Hence $j$ is Koszul-dualizing.
Since $\Phi$ is essentially surjective, Proposition \ref{prop:reflexivity} now implies that $\mathscr A$ is reflexive.
\end{proof}

The following reflexivity criteria supply the two cases in Corollary \ref{cor:main}.

\begin{proposition}[{\cite[Corollaries 3.4.3, 3.3.2, 3.3.6]{Booth-Goodbody-Opper2025}}]\label{prop:algebraic-reflexivity}
A dg algebra $A$, or a strictly unital $A_\infty$ algebra through its dg model, is reflexive in either of the following cases:
\begin{enumerate}
    \item[(i)] $H^{>0}(A)=0$ and $\dim_{\Bbbk}H^d(A)<\infty$ for every $d$;
    \item[(ii)] $H^{<0}(A)=0$, $H^0(A)$ is semisimple, and either $A$ is proper, or $A$ is smooth with $\dim_{\Bbbk}H^d(A)<\infty$ for every $d$.
\end{enumerate}
\end{proposition}

Now we are ready to prove Corollary \ref{cor:main}.

\begin{proof}[Proof of Corollary \ref{cor:main}]
Apply Proposition \ref{prop:algebraic-comparison}, Proposition \ref{prop:algebraic-reflexivity}, Morita invariance \cite[Proposition 3.8]{Goodbody2026} and Theorem \ref{thm:main}.
\end{proof}

To check the generation condition in \eqref{eq:algebraic-density-test}, let $R=\Bbbk^r$, with standard idempotents $e_1,\ldots,e_r$, and put $R_i=e_iR$.
A strictly unital minimal $A_\infty$ algebra $B$ over $R$ is \emph{positive} if $B^{<0}=0$ and $B^0=R$; projection to $R$ is an augmentation \cite[Section 4.2]{Su-Yang2019}.

\begin{proposition}\label{prop:algebraic-generation}
One has $\mathsf{PsPerf}(B)=\mathrm{thick}_B(R)$ in either of the following cases:
\begin{enumerate}
    \item[(i)] $B$ is augmented over $R$, $H^{>0}(B)=0$, $H^0(B)$ is finite-dimensional, and its augmentation identifies $H^0(B)/\operatorname{rad}H^0(B)$ with $R$;
    \item[(ii)] $B$ is positive.
\end{enumerate}
In either case, a module whose cohomology is $R_i[-d]$ as a graded $R$-module is quasi-isomorphic to the corresponding shifted augmentation module.
\end{proposition}
\begin{proof}
Case (i) follows by truncation and composition series over $H^0(B)$ \cite[Proposition 2.2.2 and Corollary 2.2.3]{Booth-Goodbody-Opper2025}.
Case (ii) and the assertion about module structures are \cite[Lemmas 4.2--4.3]{Su-Yang2019}.
\end{proof}

Consequently, once \eqref{eq:algebraic-action} is a quasi-isomorphism, essential surjectivity is equivalent to $R\in\mathrm{thick}_B(M)$ in these cases.
For example, suppose $G=\bigoplus_jG_j$ and $K=\bigoplus_iK_i$, with $e_i$ corresponding to $K_i$.
Writing $M_j=\mathrm{Hom}_{\mathscr B}(K,\Phi(G_j))$, formula \eqref{eq:Phi-values} gives
\[
M_je_i\simeq\mathrm{Hom}_{\mathscr A}(G_j,j_!K_i)^\vee.
\]
If $H^*(M_j)$ is $R_{i_j}[-d_j]$ as a graded $R$-module and every index $1,\ldots,r$ occurs among the $i_j$, Proposition \ref{prop:algebraic-generation} gives
\begin{equation}\label{eq:algebraic-simple-images}
M\simeq\bigoplus_jR_{i_j}[-d_j],\qquad
\mathrm{thick}_B(M)=\mathsf{PsPerf}(B).
\end{equation}

\subsection{Koszul duals of proper algebras}

The following proposition identifies the subcategory that must be divided out when starting from the other direction of the Koszul duality functor.

\begin{proposition}\label{prop:algebraic-boundary}
Let $B$ be proper and let $P$ split-generate $\mathsf{PsPerf}(B)$.
Suppose that $A=\mathrm{Hom}_B(P,P)$ is smooth and $\dim_{\Bbbk}H^d(A)<\infty$ for every $d$.
Then $A$ is reflexive, and the quasi-equivalence $\Psi=\mathrm{Hom}_B(P,-)$ identifies the following pairs of dg categories:
\begin{equation}\label{eq:algebraic-pairs}
\Psi: \bigl(\mathsf{PsPerf}(B),\mathrm{thick}_B(B^\vee)\bigr)
\xrightarrow{\;\sim\;}
\bigl(\mathsf{Perf}(A),\mathsf{PsPerf}(A)\bigr).
\end{equation}
Its inverse is $-\otimes_A P$, and it induces quasi-equivalences
\begin{equation}\label{eq:proper-algebra-boundary}
\widehat{\mathsf{Perf}(A)}_\infty
\simeq\mathsf{PsPerf}(B)/\mathrm{thick}_B(B^\vee)
\simeq\mathscr D_{\mathsf{sg}}(B^{\mathrm{op}})^{\mathrm{op}}.
\end{equation}
If $B$ satisfies the Gorenstein condition \eqref{eq:Gorenstein-cond}, the first quotient is identified with $\mathscr D_{\mathsf{sg}}(B)$.
\end{proposition}

\begin{proof}
Since $P$ split-generates $\mathsf{PsPerf}(B)$, the functor $\Psi$ is a quasi-equivalence with quasi-inverse $-\otimes_A P$.
Tensor--Hom adjunction gives
\[
\Psi(B^\vee)
=\mathrm{Hom}_B(P,B^\vee)
\simeq P^\vee.
\]
Consequently, the identification of subcategories in \eqref{eq:algebraic-pairs} reduces to proving
\[
\mathsf{PsPerf}(A)=\mathrm{thick}_A(P^\vee).
\]
We establish this equality and reflexivity of $A$ by verifying the hypotheses of Proposition \ref{prop:algebraic-comparison}.

Since $B$ is proper, $B^\vee$ belongs to $\mathsf{PsPerf}(B)$, so the preceding identification shows that $P^\vee$ is perfect as a right $A$-module.
The $A$--$B$-bimodule structure on $P$ gives $P^\vee$ a $B$--$A$-bimodule structure.
Under $\Psi(B^\vee)\simeq P^\vee$, its left $B$-action corresponds to the left action on $B^\vee$, and its right $A$-action corresponds to precomposition on $\mathrm{Hom}_B(P,B^\vee)$.
Linear duality and quasi-full faithfulness of $\Psi$ therefore give quasi-isomorphisms of dg algebras
\[
B\xrightarrow{\;\sim\;}
\mathrm{Hom}_B(B^\vee,B^\vee)
\xrightarrow{\;\sim\;}
\mathrm{Hom}_A(P^\vee,P^\vee).
\]
Thus the generators $A$ of $\mathsf{Perf}(A)$ and $P^\vee$ of $\mathrm{thick}_A(P^\vee)$ have endomorphism algebras $A$ and $B$, respectively.
These are the algebra models used in Proposition \ref{prop:algebraic-comparison}, and the corresponding graph bimodule is
\[
\mathrm{Hom}_A(A,P^\vee)\simeq {}_B P^\vee_A.
\]
Since $P \in \mathsf{PsPerf}(B)$, ${}_B P^\vee_A$ is also perfect over $\Bbbk$.
This verifies the pseudo-perfectness hypothesis on the graph bimodule.

The dual graph bimodule is $(P^\vee)^\vee\simeq P$ as an $A$--$B$-bimodule.
Hence the comparison functor in Proposition \ref{prop:algebraic-comparison} is $-\otimes_A P$, the quasi-inverse of $\Psi$ already established above.
Its action morphism \eqref{eq:algebraic-action} is the defining identification $A=\mathrm{Hom}_B(P,P)$, and the density condition \eqref{eq:algebraic-density-test} is the assumption $\mathsf{PsPerf}(B)=\mathrm{thick}_B(P)$.
Together with the assumed smoothness and degreewise finiteness of $A$, these verifications allow us to apply Proposition \ref{prop:algebraic-comparison}.
It proves that $A$ is reflexive and supplies the Koszul duality hypotheses of Lemma \ref{lem:Koszul-generation}.
The generation statement \eqref{eq:Koszul-generation} then gives
\[
\mathsf{PsPerf}(A)
=\mathrm{thick}_A(P^\vee)
=\Psi\bigl(\mathrm{thick}_B(B^\vee)\bigr),
\]
which proves \eqref{eq:algebraic-pairs}.

Passing to dg quotients and applying \eqref{eq:perf-boundary-quotient} gives the first quasi-equivalence in \eqref{eq:proper-algebra-boundary}.
Linear duality sends $\mathrm{thick}_B(B^\vee)$ to $\mathsf{Perf}(B^{\mathrm{op}})^{\mathrm{op}}$, so Lemma \ref{lem:Dsg-op-op} gives the second.
Finally, under the Gorenstein condition \eqref{eq:Gorenstein-cond}, one has $\mathrm{thick}_B(B^\vee)=\mathsf{Perf}(B)$, and the first quotient is therefore $\mathscr D_{\mathsf{sg}}(B)$.
\end{proof}

\begin{corollary}\label{cor:positive-smooth-dual}
Let $B$ be a finite-dimensional positive $A_\infty$ algebra over $R$ with $B^1=0$.
Then $A=\mathrm{Hom}_B(R,R)$ has a nonpositive smooth dg model with $H^0(A)=R$ and finite-dimensional cohomology in every degree.
The conclusions of Proposition \ref{prop:algebraic-boundary} hold with $P=R$.
\end{corollary}
\begin{proof}
Put $V=(B^{>0}[1])^\vee$.
The relative bar construction gives a dg model for $A$ whose underlying graded algebra is $T_R(V)$; compare \cite[Section 2.2 and Lemma 5.3(b)]{Lu-Palmieri-Wu-Zhang2008}.
Since $B$ is finite-dimensional and $B^1=0$, the space $V$ is finite-dimensional and concentrated in strictly negative degrees.
Consequently, only finitely many tensor lengths contribute to each graded piece. In particular, this model is nonpositive and has finite-dimensional cohomology in every degree.

The augmentation excludes constant terms, and minimality of $B$ excludes linear terms in the differential of each generator.
Thus, for a homogeneous generator $v\in V$, the differential $d(v)$ is a finite sum of words $v_1\cdots v_\ell$ with $\ell\geqslant2$.
Each such word has degree $|v|+1$, and all its factors have negative degree.
It follows that
\[
|v_i|
=|v|+1-\sum_{j\ne i}|v_j|
>|v|
\qquad\text{for every }i.
\]
Hence $d(v)$ involves only generators of strictly greater degree than $v$.

To deduce smoothness of $A$, let $\delta:A\to\Omega^1_{A/R}$ be the universal relative derivation.
As a graded $A$-bimodule, one has
\[
\Omega^1_{A/R}\cong A\otimes_R V\otimes_R A,
\]
with $1\otimes v\otimes1$ corresponding to $\delta v$.
Its differential satisfies
\[
d(\delta v)=\delta(d(v)).
\]
By the Leibniz rule, the right-hand side is a sum of bimodule multiples of elements $\delta v_i$ with $|v_i|>|v|$.
Thus, when the generators $\delta v$ are adjoined in decreasing order of degree, the differential of each new generator lies in the bimodule generated by those already adjoined.
Since $V$ is finite-dimensional, this gives a finite filtration of $\Omega^1_{A/R}$ by dg sub-bimodules whose successive quotients are finite direct sums of shifts of $Ae_i\otimes_{\Bbbk}e_jA$.
These quotients are perfect, so $\Omega^1_{A/R}$ is a perfect $A$-bimodule.

Since $R=\Bbbk^r$, the bimodule $A\otimes_R A$ is also a finite direct sum of perfect bimodules $Ae_i\otimes_{\Bbbk}e_iA$.
The exact sequence
\[
0\longrightarrow\Omega^1_{A/R}
\longrightarrow A\otimes_R A
\longrightarrow A
\longrightarrow0
\]
therefore shows that the diagonal bimodule $A$ is perfect, hence smooth.
Moreover, $A^0=R$, and $d(A^{-1})=0$ because the differential of a degree $-1$ generator has no constant term.
Thus $H^0(A)=R$.

Finally, Proposition \ref{prop:algebraic-generation}(ii) gives $\mathsf{PsPerf}(B)=\mathrm{thick}_B(R)$.
Since $B$ is proper and $A$ is smooth with finite-dimensional cohomology in every degree, Proposition \ref{prop:algebraic-boundary} applies with $P=R$.
\end{proof}

\begin{remark}\label{rem:bar-finiteness}
Degreewise finiteness can also hold when $B^1\ne0$.
For a finite-dimensional positive $B$, the reduced bar construction is
\[
\mathrm{Bar}(B)=\bigoplus_{\ell\geqslant0}(B^{>0}[1])^{\otimes_R\ell},\qquad
\mathrm{Hom}_B(R,R)\simeq(\mathrm{Bar}(B)^\vee)^{\mathrm{op}},
\]
where the opposite is taken for the forward convolution product.
If the degree-one arrows, indexed by bases of $e_jB^1e_i$ from $i$ to $j$, form an acyclic quiver, blocks of degree-zero bar letters have bounded length; see \cite[Section II.1 and Lemma II.1.4(c)]{Assem-Simson-Skowronski2006}.
Each graded piece of $\mathrm{Bar}(B)$ is then finite-dimensional, giving the finiteness required in Proposition \ref{prop:algebraic-comparison}.
\end{remark}

\subsection{Geometric examples}\label{subsec:geometric-examples}

For a finite-type Weinstein manifold $X^{2n}$, let $\mathscr W(X)$ and $\mathscr F(X)$ denote its wrapped and compact Fukaya categories of exact graded branes with finite-rank local systems.
$\mathscr W(X)$ is smooth with a smooth $n$-Calabi-Yau structure and is generated by cocores, while Poincar\'{e} duality gives a proper $n$-Calabi--Yau structure on $\mathscr F(X)$; see \cite{Ganatra2013,Ganatra-Pardon-Shende2024a,Chantraine-DimitroglouRizell-Ghiggini-Golovko2024} and \cite[Theorem 2]{Ganatra2023}.

The following examples illustrate the algebraic criteria and the scheme comparison \eqref{eq:scheme-boundary-comparison}, giving
\begin{equation}\label{eq:geometric-specialization}
\widehat{\mathsf{Perf}(\mathscr W(X))}_\infty
\simeq\mathscr D_{\mathsf{sg}}(\mathscr F(X)).
\end{equation}
By \cite[Theorem 1.1]{Ganatra-Gao-Venkatesh2022}, this also identifies the Rabinowitz Fukaya category up to Morita equivalence.
In the first two examples, we first apply the criteria to the full subcategory on the chosen compact cores and then verify that these cores split-generate $\mathsf{Perf}(\mathscr F(X))$.

\begin{enumerate}
    \item \emph{Cotangent bundles.}
    Let $Q$ be a closed simply connected spin manifold, with coefficients in any field and zero background class.
    A cotangent fibre $G$ and the zero section $K$ give
    \[
    A\simeq C_{-*}(\Omega Q;\Bbbk),\qquad
    B\simeq C^*(Q;\Bbbk),\qquad M\simeq\Bbbk;
    \]
    see \cite[Theorem 1.1]{Abouzaid2011} and \cite[Introduction]{Ekholm-Lekili2023} for the compatible Yoneda and Koszul actions.
    A finite minimal model of $B$ is positive with $B^1=0$, so Corollary \ref{cor:positive-smooth-dual} applies.
    Since $A$ is nonpositive with $H^0(A)=\Bbbk$, Proposition \ref{prop:algebraic-generation}(i) places every compact object's fibre Yoneda module in $\mathrm{thick}_A(\Bbbk)$.
    Thus the zero section split-generates $\mathsf{Perf}(\mathscr F(T^*Q))$, giving \eqref{eq:geometric-specialization} for the full compact category.

    \item \emph{Tree plumbings.}
    Let $X$ be a finite tree plumbing of cotangent bundles of closed simply connected spin $n$-manifolds, $n\geqslant3$, with the same coefficient convention.
    Take $G$ to be the sum of the cocores and $K$ the sum of the cores.
    Grade the cores so that opposite intersection generators have degrees $1,n-1$.
    Then $B$ is proper and admits a positive minimal model with acyclic degree-one quiver, and the summands of $M$ are the augmentation modules $R_i$.
    The action map \eqref{eq:algebraic-action} is a quasi-isomorphism by \cite[Theorems 66, 68, and 83; Section 6.1.4]{Ekholm-Lekili2023}.
    Proposition \ref{prop:algebraic-generation} and Remark \ref{rem:bar-finiteness} verify \eqref{eq:algebraic-density-test}.
    Moreover, $A$ is nonpositive with finite-dimensional $H^0(A)$, and the cocore Yoneda modules of the cores represent all simple $H^0(A)$-modules.
    The same truncation argument therefore proves compact generation and yields \eqref{eq:geometric-specialization}.

    \item \emph{Punctured tori.}
    Over a perfect field, the standard graded torus $T_r$ with $r\geqslant1$ punctures has a quasi-equivalence of pairs
    \[
    \bigl(\mathsf{Perf}(\mathscr W(T_r)),\mathsf{Perf}(\mathscr F(T_r))\bigr)
    \simeq
    \bigl(\mathsf{Coh}(C_r),\mathsf{Perf}(C_r)\bigr),
    \]
    where $C_1$ is the split nodal cubic and $C_r$ is the cycle of $r$ projective lines for $r>1$; see \cite[Theorem 9.1]{Lekili-Perutz2015} and \cite[Corollary 3.4.1 and Theorem 3.7.2]{Lekili-Polishuk2017}.
    Since its dualizing complex is $\mathcal O_{C_r}[1]$, formula \eqref{eq:scheme-boundary-comparison} gives \eqref{eq:geometric-specialization}, with target $\mathscr D_{\mathsf{sg}}(C_r)$.
\end{enumerate}

Calculations for double-bubble plumbings appear in \cite[Proposition 5.2 and Corollaries 5.7--5.8]{Smith-Wemyss2023}, and for Milnor fibres in \cite[Section 3]{Lekili-Ueda2021} and \cite[Theorem 6.11]{Lekili-Ueda2022}.
For the description of double suspensions by cluster categories, see \cite[Theorems 2.2 and 4.1]{Lekili-Ueda2024}.

\section{Hochschild cohomology}\label{sec:hochschild}

We deduce multiplicative comparisons of Hochschild cohomology from Theorem \ref{thm:main}, then relate them to singular Hochschild cohomology under Keller's hypotheses.
All module Hom and tensor products use the dg models fixed in Section \ref{sec:basic}.

\subsection{The boundary comparison}

For a small dg category $\mathscr C$ and a $\mathscr C$-bimodule $M$, write
\[
C^*(\mathscr C;M)
=\mathrm{Hom}_{\mathscr C^{\mathrm e}}(\mathscr C_\Delta,M),
\qquad
HH^m(\mathscr C,M)=H^mC^*(\mathscr C;M),
\]
and omit $M$ for diagonal coefficients; see \cite[Section 5.4]{Keller2006b} for the Hochschild bar model and cup product.
Let $q:\mathscr C\to\mathscr D$ be a dg localization, meaning that restriction $q^*:\mathsf{Mod}(\mathscr D)\to\mathsf{Mod}(\mathscr C)$ is quasi-fully faithful; dg quotients have this property \cite[Section 3.9]{Keller2011}.
Set
\begin{equation}\label{eq:boundary-coefficient}
\partial_q\mathscr C=(q^{\mathrm e})^*\mathscr D_\Delta,
\qquad
(\partial_q\mathscr C)(X,Y)=\mathscr D(qY,qX),
\qquad
q^{\mathrm e}=q^{\mathrm{op}}\otimes q.
\end{equation}
Composition in $\mathscr D$ and the map $\eta:\mathscr C_\Delta\to\partial_q\mathscr C$ induced by $q$ give this bimodule a multiplication and unit, hence a cup product on $HH^*(\mathscr C,\partial_q\mathscr C)$.

\begin{proposition}\label{prop:HH-localization}
Restriction gives a quasi-isomorphism
\begin{equation}\label{eq:HH-localization}
q^*:C^*(\mathscr D)\longrightarrow C^*(\mathscr C;\partial_q\mathscr C)
\end{equation}
and an isomorphism of unital graded algebras on cohomology.
\end{proposition}
\begin{proof}
By \cite[Proposition 3.10(a)]{Keller2011}, extension of scalars along $q^{\mathrm e}$ sends $\mathscr C_\Delta$ to $\mathscr D_\Delta$.
Adjunction therefore gives the quasi-isomorphism \eqref{eq:HH-localization}.
On bar cochains, its arity-$m$ component is
\[
(q^*f)_r(a_m,\ldots,a_1)
=f_m(q(a_m),\ldots,q(a_1))
\]
for composable morphisms $a_1,\ldots,a_m$ in $\mathscr C$.
The cup product uses deconcatenation of these inputs and composition in $\mathscr D$, so this formula preserves products and sends the diagonal unit to $\eta$.
\end{proof}

Under the hypotheses of Theorem \ref{thm:main}, apply this proposition to
\[
q:\mathscr C=\mathsf{Perf}(\mathscr A)
\longrightarrow\mathscr D
=\mathsf{Perf}(\mathscr A)/\mathsf{PsPerf}(\mathscr A).
\]
For the Yoneda functor $i:\mathscr A\to\mathscr C$, define the \emph{boundary coefficient bimodule} by
\begin{equation}\label{eq:boundary-coefficient-A}
\partial_\infty\mathscr A=(i^{\mathrm e})^*\partial_q\mathscr C.
\end{equation}
Proposition \ref{prop:boundary-quotient} identifies it with the pullback of the diagonal of $\widehat{\mathscr A}_\infty$ along \eqref{eq:reduced-Yoneda}.

\begin{corollary}\label{cor:HH-boundary}
Under the hypotheses of Theorem \ref{thm:main}, there are natural isomorphisms of unital graded algebras
\begin{equation}\label{eq:HH-at-infinity}
HH^*(\mathscr A,\partial_\infty\mathscr A)
\cong HH^*(\widehat{\mathscr A}_\infty)
\cong HH^*(\widehat{\mathsf{Perf}(\mathscr A)}_\infty).
\end{equation}
If $\mathscr A$ is reflexive, these algebras are isomorphic to $HH^*(\mathscr D_{\mathsf{sg}}(\mathscr B^{\mathrm{op}})^{\mathrm{op}})$.
If $\mathscr A$ is reflexive and $\mathscr B$ satisfies the Gorenstein condition \eqref{eq:Gorenstein-cond}, the functor \eqref{eq:ordinary-comparison} induces an isomorphism of unital graded algebras
\begin{equation}\label{eq:HH-singularity}
HH^*(\mathscr D_{\mathsf{sg}}(\mathscr B))
\xrightarrow{\;\sim\;}HH^*(\mathscr A,\partial_\infty\mathscr A).
\end{equation}
\end{corollary}
\begin{proof}
The Yoneda Morita equivalence preserves Hochschild cohomology with corresponding coefficients and their cup products by \cite[Proposition 3.10(e)]{Keller2011}.
Together with Proposition \ref{prop:HH-localization}, it identifies $HH^*(\mathscr A,\partial_\infty\mathscr A)$ with $HH^*(\mathscr D)$.
Proposition \ref{prop:boundary-quotient} gives \eqref{eq:HH-at-infinity}, and Theorem \ref{thm:main} gives the comparisons with the singularity categories.
\end{proof}

After translating the opposite-category and grading conventions, this is the cup product of \cite[Definition 4.2 and Remark 4.3]{Rivera-Takeda-Wang2024}; the boundary restriction map is the unital map of \cite[Proposition 3.10]{Gao2025}.
The same argument applies to the algebraic pairs \eqref{eq:algebraic-pairs} and the geometric equivalences \eqref{eq:geometric-specialization}.

\subsection{Singular Hochschild cohomology comparison}

For proper $\mathscr B$, its diagonal is pseudo-perfect over $\mathscr B^{\mathrm e}$ by Definition \ref{def:smooth-proper}.
Denote its image in $\mathscr D_{\mathsf{sg}}(\mathscr B^{\mathrm e})$ by $\overline{\mathscr B}_\Delta$, and define
\[
HH^m_{\mathsf{sg}}(\mathscr B)
=H^m\mathrm{Hom}_{\mathscr D_{\mathsf{sg}}(\mathscr B^{\mathrm e})}
(\overline{\mathscr B}_\Delta,\overline{\mathscr B}_\Delta),
\]
with product induced by composition; for ordinary finite-dimensional algebras, this is singular Hochschild cohomology as in \cite[Section 1]{Keller2020}.

For a finite-dimensional ordinary algebra $B$ with $B/\operatorname{rad}(B)$ separable over $\Bbbk$, Keller's canonical comparison is an isomorphism of unital graded algebras
\begin{equation}\label{eq:singular-HH-comparison}
HH^*_{\mathsf{sg}}(B)\xrightarrow{\;\sim\;}HH^*(\mathscr D_{\mathsf{sg}}(B)).
\end{equation}
Indeed, \cite[Theorem 3.7]{Elagin-Lunts-Schnurer2020} gives smoothness of $\mathsf{PsPerf}(B)$, which is the hypothesis required by \cite[Theorem 1.1]{Keller2020}.
Consequently, when $\mathscr B=B$ in Corollary \ref{cor:HH-boundary}, $\mathscr A$ is reflexive, and $B$ satisfies \eqref{eq:Gorenstein-cond}, composing \eqref{eq:singular-HH-comparison} with \eqref{eq:HH-singularity} identifies $HH^*_{\mathsf{sg}}(B)$ with $HH^*(\mathscr A,\partial_\infty\mathscr A)$ as unital graded algebras.
If $B$ is split, meaning $B/\operatorname{rad}(B)\cong\prod_i M_{d_i}(\Bbbk)$, \cite[Theorem 1.1]{Xi-Zhang2026}, applied to $B^{\mathrm{op}}$ for our right-module convention, lifts \eqref{eq:singular-HH-comparison} to the homotopy category of $B_\infty$-algebras.

The complete multiplicative solution to Conjecture \ref{conjecture:Takeda} requires further structure on $\mathscr A$.
For example, $\mathscr A$ also carries a nondegenerate $4$-truncated $n$-pre-Calabi--Yau structure, in the sense of \cite[Definitions 2.9 and 2.11]{Rivera-Takeda-Wang2024}, its induced Calabi--Yau duality identifies
\[
HH^r(\mathscr A,\partial_\infty\mathscr A)
\cong HH_{n-r}(\mathscr A,\partial_\infty\mathscr A).
\]
With the Koszul signs for this regrading, \cite[Proposition 5.8 and Theorem 5.11]{Rivera-Takeda-Wang2024} identify the cup product with their product of homological degree $-n$.

\printbibliography

\end{document}